\documentclass[12pt]{article}

\usepackage[margin=1in]{geometry}
\usepackage{graphicx}
\usepackage{amsmath}
\usepackage{amsfonts}
\usepackage{amsthm}
\usepackage{amssymb,latexsym,hyperref}

\newtheorem{thm}{Theorem}[section]
\newtheorem{lem}[thm]{Lemma}
\newtheorem{claim}[thm]{Claim}
\newtheorem{defn}[thm]{Definition}
\newtheorem{ques}[thm]{Question}

\newtheorem{rem}[thm]{Remark}

\def\diam{\diamondsuit}
\def\club{\clubsuit}

\def\res{\upharpoonright}
\def\name{\mathring}

\def\dom{\textsf{dom}}

\def\bb{\mathbb}
\def\cal{\mathcal}
\def\frak{\mathfrak}
\def\e{\varepsilon}
\def\club{\clubsuit}
\def\diam{\diamondsuit}

\def\cf{\textsf{cf}}
\def\dom{\textsf{dom}}

\def\Lim{\textsf{•}{Lim}}
\def\otp{\textsf{otp}}

\def\max{\text{max}}
\def\min{\text{min}}
\def\sf{\textsf}
\def\CR{\sf{CR}}
\def\pos{\sf{pos}}

\def\diam{\diamondsuit}

\begin{document}

\nocite{*}

\title{On some variants of the club principle}

\author{Ashutosh Kumar\footnote{Einstein Institute of Mathematics, The Hebrew University of Jerusalem, Edmond J Safra Campus, Givat Ram, Jerusalem 91904, Israel; email: \textsl{akumar@math.huji.ac.il}; Supported by a Postdoctoral Fellowship at the Einstein Insititute of Mathematics funded by European Research Council grant 338821}, Saharon Shelah\footnote{Einstein Institute of Mathematics, The Hebrew University of Jerusalem, Edmond J Safra Campus, Givat Ram, Jerusalem 91904, Israel and Department of Mathematics, Rutgers, The State University of New Jersey, Hill Center-Busch Campus, 110 Frelinghuysen Road, Piscataway, NJ 08854-8019, USA; email: \textsl{shelah@math.huji.ac.il}; Partially supported by European Research Council grant 338821; Publication no. 1136}}

\date{}

\maketitle

\begin{abstract}

We study some asymptotic variants of the club principle. Along the way, we construct some forcings and use them to separate several of these principles.

\end{abstract}

\section{Introduction}

For a regular uncountable cardinal $\kappa$ and a stationary $S \subseteq \Lim(\kappa)$, the club principle $\club_S$ says the following: There exists $\bar{A} = \langle A_{\delta} : \delta \in S\rangle$ where each $A_{\delta}$ is an unbounded subset of $\delta$ of order type $\cf(\delta)$ such that for every $A \in [\kappa]^{\kappa}$, there exists some (equivalently, stationary many) $\delta \in S$ such that $A_{\delta} \subseteq A$. We say that $\bar{A}$ is a $\club_S$ witnessing sequence. If $\kappa = \omega_1$ and $S = \Lim(\omega_1)$ is the set of all countable limit ordinals, we drop the $S$ and write $\club$. \\

In $\cite{DjSh574}$, it was shown that $\club^1$ does not imply $\club$ where $\club^1$ is the following statement: There exists $\bar{A} = \langle A_{\delta} : \delta \in \Lim(\omega_1) \rangle$ where each $A_{\delta}$ is an unbounded subset of $\delta$ of order type $\omega$ such that for every $A \in [\omega_1]^{\aleph_1}$, there exists $\delta$ such that $A_{\delta} \setminus A$ is finite. For some other variants of the club principle, see \cite{DjSh574, KmSh1046, KmSh1063}. \\

\begin{defn}
For $a \in (0, 1]$ and a stationary set $S \subseteq \Lim(\omega_1)$, the principle $\club^{\inf \geq a}_S$ says the following: There exists $\bar{A} = \langle A_{\delta} : \delta \in S \rangle$ such that 

\begin{itemize}

\item[(a)] each $A_{\delta} = \{\alpha_{\delta, n} : n < \omega \}$ and $\alpha_{\delta, n}$'s are increasing cofinal in $\delta$ and

\item[(b)] for every $A \in [\omega_1]^{\aleph_1}$, there exists $\delta \in S$ such that $$ \liminf_n \frac{|\{k < n: \alpha_{\delta, k} \in A\}|}{n} \geq a$$

\end{itemize}

If $S = \Lim(\omega_1)$, we write $\club^{\inf \geq a}$. By $\club^{\lim}$, we mean $\club^{\inf \geq 1}$. 
\end{defn}

It is clear that $\club^1$ implies $\club^{\lim}$ and for $0 < a < b \leq 1$, $\club^{\inf \geq b}$ implies $\club^{\inf \geq a}$. At the end of Section 1, we show that under CH, all of these principles are equivalent to diamond. 

\begin{thm}
\label{clubch}
Assume CH. Then for every $a \in (0, 1]$, $\club^{\inf \geq a}$ implies $\diam$.
\end{thm}

The bulk of the work in this paper is to show the following.

\begin{thm} 
\label{infthms}
\item[(1)] $\club^{\lim} \wedge \neg \club^1$ is consistent.
\item[(2)] For every $a \in (0, 1]$, $\club^{\inf \geq a} \wedge (\forall b > a) \neg \club^{\inf \geq b}$ is consistent.
\item[(3)] For every $a \in (0, 1]$, $\neg \club^{\inf \geq a} \wedge (\forall b < a) \club^{\inf \geq b}$ is consistent.

\end{thm}

In Sections 2-5 we prove Theorem \ref{infthms}(1). In Section 6, we supply the necessary modifications to get parts (2) and (3). The forcing used is quite flexible and can be useful for separating many similar principles. \\

In Section 7, we introduce $\club^{\sup \geq a}$ (defined analogously) and prove the following in ZFC.

\begin{thm}
\label{supequiv}
For every $a, b \in (0, 1)$, $\club^{\sup \geq a}$ is equivalent to $\club^{\sup \geq b}$.
\end{thm}

Finally, in Section 8, we prove that

\begin{thm}
\label{supsep}
$\club^{\sup \geq 0.5} \wedge \neg \club^{\sup \geq 1}$ is consistent.
\end{thm}

\textbf{On notation}: $\Lim(\kappa)$ denotes the set of all limit ordinals below $\kappa$. $\cf(\alpha)$ is the cofinality of $\alpha$. $S^{\kappa}_{\delta} = \{\alpha < \kappa : \cf(\alpha) = \cf(\delta)\}$. For $k \leq \omega$, $\omega^k$ is the $k$th ordinal power of $\omega$ with under ordinal exponentiation. For $a, b$ sets of ordinals, then we write $a < b$ to denote $(\forall \alpha \in a)(\forall \beta \in b)(\alpha < \beta)$. In forcing, we use the convention that a larger condition is the stronger one - $p \geq q$ means $p$ extends $q$.

\subsection{CH and $\club^{\inf}$}

Recall that $\diam$ says the following: There exists $\langle A_{\delta}: \delta \in \Lim(\omega_1) \rangle$ where each $A_{\delta} \subseteq \delta$ such that for every $A \in [\omega_1]^{\aleph_1}$, $\{\delta \in \Lim(\omega_1): A_{\delta} = A \cap \delta \}$ is stationary. An equivalent formulation (see \cite{Kun}) is the following: There exists $\langle \cal{A}_{\delta}: \delta \in \Lim(\omega_1) \rangle$ where each $\cal{A}_{\delta}$ is a countable family of subsets of $\delta$ such that for every $A \in [\omega_1]^{\aleph_1}$, $\{\delta \in \Lim(\omega_1):  A \cap \delta \in \cal{A}_{\delta} \}$ is stationary. \\

Proof of Theorem \ref{clubch}: Assume CH. Suppose $a \in (0, 1]$ and $\club^{\inf \geq a}$ holds as witnessed by $\bar{A} = \langle A_{\delta}: \delta \in \Lim(\omega_1) \rangle$. Let $A_{\delta} = \{\alpha_{\delta, n}: n < \omega_1\}$ list $A_{\delta}$ in increasing order. Using CH, fix $\langle B_i: i < \omega_1 \rangle$ such that each $B_i \subseteq i$ and for every $B \in [\omega_1]^{\leq \aleph_0}$, there are uncountably many $i < \omega_1$ for which $B = B_i$. \\

For $\delta \in \Lim(\omega_1)$, define $\cal{A}_{\delta}$ as follows. $A \in \cal{A}_{\delta}$ iff for some $u \subseteq \omega$ the following hold.

\begin{itemize}

\item[(a)] $\liminf_n |u \cap n|/n \geq a$.

\item[(b)] For every $m < n$ in $u$, $B_{\alpha_{\delta, m}} = B_{\alpha_{\delta, n}} \cap \alpha_{\delta, m}$ and $A = \bigcup_{n \in u} B_{\alpha_{\delta, n}} $.

\end{itemize}

We claim that each $\cal{A}_{\delta}$ is finite. In fact, $|\cal{A}_{\delta}| \leq 1/a$. To see this assume otherwise and let $\{A_k : k < K\}$ be pairwise distinct members of $\cal{A}_{\delta}$ where $Ka > 1$. Choose $\langle u_k: k < K \rangle$ witnessing $A_k \in \cal{A}_k$. Choose $N_1 < N_2$ such that the following hold.

\begin{itemize}

\item[(i)] $\langle A_k \cap \alpha_{\delta, N_1}: k < K \rangle$ has pairwise distinct members

\item[(ii)] $|u_k \cap [N_1, N_2)| > (N_2 - N_1)/K$ for each $k < K$

\end{itemize}

By (ii), it follows that for some $j < k < K$, $[N_1, N_2) \cap u_j \cap u_k \neq \emptyset$. But if $n \in [N_1, N_2) \cap u_j \cap u_k$, then $B_{\alpha_{\delta, n}} = A_j \cap \alpha_{\delta, n} = A_k \cap \alpha_{\delta, n}$ which is impossible by (i). \\

To complete the proof it is enough to show the following.

\begin{claim}
\label{diam}
For every $X \in [\omega_1]^{\aleph_1}$, for every club $E \subseteq \omega_1$, there exists $\delta \in E$ such that $C \cap \delta \in \cal{A}_{\delta}$.
\end{claim}

Proof of Claim \ref{diam}: Construct $\langle \alpha_i : i < \omega_1 \rangle$ such that $\alpha_i$'s are increasing and for every $i < \omega_1$, $X \cap \sup_{j < i} \alpha_j = B_{\alpha_i}$. Choose $\delta \in E$ and $u \subseteq \omega$ such that $\liminf_n |u \cap n|/n \geq a$ and $\{\alpha_{\delta, n}: n \in u\} \subseteq \{\alpha_i : i < \omega_1 \}$. It follows that $X \cap \delta = \bigcup_{n \in u} B_{\alpha_{\delta, n}} \in \cal{A}_{\delta}$. \qed

\section{Creatures}
\label{c1}
Fix a family $\{S_k: k < \omega\}$ of pairwise disjoint stationary subsets of $\omega_1$ consisting of limit ordinals. We describe a ccc forcing which is somewhat intermediate between adding $\aleph_1$ Cohen reals and adding a Cohen subset of $\omega_1$.

\begin{defn} 
\label{crdefn}
We say that $(\CR, \Sigma)$ is an $\aleph_1$-CP (creating pair) if the following hold.

\begin{itemize}

\item[(A)] We call members of $\CR$ creatures. For each $\frak{c} \in \CR$, 

\subitem (i) $\frak{c} = (\dom(\frak{c}), \pos(\frak{c}), f_{\frak{c}})$.

\subitem (ii) $\dom(\frak{c})$ is a non empty subset of $\omega_1$ of order type $< \omega^{\omega}$.

\subitem (iii) For every limit $\delta < \omega_1$, if $\dom(\frak{c}) \cap \delta$ is unbounded in $\delta$, then for some $k \geq 1$, $\delta \in S_k$ and $\otp(\dom(\frak{c}) \cap \delta) = \e + \omega^{j}$ for some $\e < \omega^{\omega}$ and $1 \leq j \leq k$ - In particular, for every $\delta \in S_0$, $\dom(\frak{c}) \cap \delta$ is bounded below $\delta$.

\subitem (iv) $\pos(\frak{c})$ (possibilities for $\frak{c}$) is a countable set of functions from $\dom(\frak{c})$ to $\{0, 1\}$ and $f_{\frak{c}} \in \pos(\frak{c})$.

\subitem (v) If $\dom(\frak{c})$ is finite, then $\pos(\frak{c}) = \{f_{\frak{c}}\}$ - We call such $\frak{c}$ finite creature.

\item[(B)] For every finite $u \subseteq \omega_1$, and $f:u \to \{0, 1\}$, there exists $\frak{c} \in \CR$ such that $\dom(\frak{c}) = u$ and $f_{\frak{c}} = f$.

\item[(C)] For every $\delta < \omega_1$, $|\{\frak{c} \in \CR : \dom(\frak{c}) \subseteq \delta\}| \leq \aleph_0$.

\item[(D)] $\Sigma$ is a function with domain $\CR$ that satisfies the following.

\subitem (i) $\Sigma(\frak{c})$ is a countable set of finite tuples $\bar{\frak{d}} = \langle \frak{d}_k : k < n \rangle$ where

\subsubitem (a) $\frak{d}_k \in \CR$,

\subsubitem (b) $\dom(\frak{c}) = \bigcup_{k < n} \dom(\frak{d}_k)$,

\subsubitem (c) $\dom(\frak{d}_k) < \dom(\frak{d}_{k+1})$ and

\subsubitem (d) whenever $f_k \in \pos(\frak{d}_k)$ for $k < n$, $\bigcup_{k < n} f_k \in \pos(\frak{c})$.

\subitem (ii) Cuts: If $\frak{c} \in \CR$ and $\alpha \in \dom(\frak{c})$ then for some $\bar{d} = \langle \frak{d}_k : k < n \rangle \in \Sigma(\frak{c})$, there exists $k < n$ such that $\min(\dom(\frak{d}_k)) = \alpha$.

\subitem (iii) $\langle \frak{c} \rangle \in \Sigma(\frak{c})$.

\subitem (iv) Transitivity: If $\langle \frak{c}_k: k < n \rangle \in \Sigma(\frak{c})$ and $\langle \frak{d}_{k, l}: l < n_k \rangle \in \Sigma(c_k)$ for $k < n$, then $\langle \frak{d}_{k, l} : k < n, l < n_k \rangle \in \Sigma(\frak{c})$.

\item[(E)] Finite joins: If $\{\frak{d}_k: k < n\} \subseteq \CR$ and $\dom(\frak{d}_k) < \dom(\frak{d}_{k+1})$, then there exists $\frak{c} \in \CR$ such that 

\subitem (i) $\dom(\frak{c}) = \bigcup_{k < n} \dom(\frak{d}_k)$,

\subitem (ii) $\pos(\frak{c}) =  \{\bigcup_{k < n} f_k: (\forall k < n)(f_k \in \pos(\frak{d}_k))\}$,

\subitem (iii) $f_{\frak{c}} = \bigcup_{k < n} f_{\frak{d}_k}$ and

\subitem (iv) $\Sigma(\frak{c}) = \{ \bigcup_{i < n} \bar{\frak{f}}_i :(\forall i < n) (\bar{\frak{f}}_i \in \Sigma(\frak{d}_i)) \} $.

\end{itemize}

\end{defn}

\begin{defn}
\label{f1}
Suppose $(\CR, \Sigma)$ is an $\aleph_1$-CP. Define $\bb{Q} = \bb{Q}_{\CR, \Sigma}$ to be the forcing whose conditions are $p = \{\frak{c}_k : k < n\}$ where $\frak{c}_k \in \CR$ and $\dom(\frak{c}_k) < \dom(\frak{c}_{k+1})$. We write $\dom(p)$ for $\bigcup_{\frak{c} \in p} \dom(\frak{c})$. For $p, q \in \bb{Q}$, define $p \leq q$ iff for every $\frak{c} \in p$, there exists $\bar{\frak{d}} = \langle \frak{d}_k : k < n \rangle \in \Sigma(\frak{c})$ such that $\{\frak{d}_k : k < n\} \subseteq q$. Define $\bb{Q} \res \alpha = \{p \in \bb{Q}: \dom(p) \subseteq \alpha\}$. Let $$\name{f}_{\bb{Q}} = \bigcup \{f_{\frak{d}}: (\exists p \in G_{\bb{Q}})(\frak{d} \in p \text{ is a finite creature})\}$$

Note that $\Vdash_{\bb{Q}} \name{f}_{\bb{Q}}: \omega_1 \to \{0, 1\}$
\end{defn}

\textbf{Example}: Let $\CR$ be the set of all finite creatures $\frak{c} = (F, \{f\}, f)$ - So $F \subseteq \omega_1$ is finite and $f: F \to \{0, 1\}$. Let $\Sigma(\frak{c})$ be the set of all $\bar{\frak{d}}$ such that the join of the members of $\bar{\frak{d}}$ is $\frak{c}$. Then forcing with $\bb{Q} = \bb{Q}_{\CR, \Sigma}$ is same as adding $\aleph_1$ Cohen reals. Note that this destroys all old witnesses to $\club^{\lim}$. We would later add more creatures to $\CR$ in such a way that while some old $\club^{\lim}$ witnessing sequences are preserved, all old $\club^1$ witnessing sequences are destroyed. \\

Recall that a forcing notion $\bb{Q}$ has $\aleph_1$ as a precaliber if whenever $\{ p_i : i < \omega_1 \} \subseteq \bb{Q}$, there exists $X \in [\omega_1]^{\aleph_1}$ such that $\{p_i: i \in X\}$ is centered - i.e., for every finite $F \subseteq X$, there exists $p \in \bb{Q}$ such that $(\forall i \in F)(p_i \leq p)$.

\begin{claim} 
\label{CRccc}
Suppose $(\CR, \Sigma)$ is an $\aleph_1$-CP. Let $\bb{Q} = \bb{Q}_{\CR, \Sigma}$. Then $\bb{Q}$ has $\aleph_1$ as a precaliber.
\end{claim}

Proof of Claim \ref{CRccc}: Suppose $\{p_i : i < \omega_1\} \in [\bb{Q}]^{\aleph_1}$. The map $i \mapsto k(i) = \sup(\bigcup_{\frak{c} \in p_i} \dom(\frak{c}) \cap i) $ is regressive on $S_0$. Choose $X_1 \in [S_0]^{\aleph_1}$ and $k(\star) < \omega_1$ such that for every $i \in X_1$, $k(i) = k(\star)$ and for every $i < j$ in $X_1$, $\dom(p_i) \cap \dom(p_j) \subseteq k(\star)$. Using Definition \ref{crdefn}(D)(ii), by possibly extending each $p_i$, we can assume that for every $\frak{c} \in p_i$, either $\dom(\frak{c}) \subseteq k(\star)$ or $\inf(\dom(\frak{c})) \geq k(\star)$. Since $\{\frak{c} \in \CR: \dom(\frak{c}) \subseteq k(\star) \}$ is countable, we can find $X \in [X_1]^{\aleph_1}$ such that for every $i \in X$, $\{\frak{c} \in p_i: \dom(\frak{c}) \subseteq k(\star)\}$ does not depend on $i \in X$. Now for any finite $F \subseteq X$, $\bigcup_{i \in F} p_i$ is a common extension of $\{p_i: i \in F\}$. \qed

\begin{claim}
\label{deltacr}

Suppose $(\CR, \Sigma)$ is an $\aleph_1$-CP. Let $\bb{Q} = \bb{Q}_{\CR, \Sigma}$. Let $\langle p_i : i < \omega_1 \rangle$ be a sequence of conditions in $\bb{Q}$ such that for every $i < j < \omega_1$, $\sup(\dom(p_i)) < \sup(\dom(p_j))$. Then there exist $X \in [\omega_1]^{\aleph_1}$, $\langle q_i : i \in X \rangle$, $m < n < \omega$ such that for every $i \in X$

\begin{itemize}

\item[(a)] $q_i \in \bb{Q}$, $q_i \geq p_i$ and $\dom(q_i) = \dom(p_i)$,

\item[(b)] $q_i = \{\frak{c}_{i, k}: k < n\}$ and for every $k < n-1$, $\dom(\frak{c}_{i, k}) < \dom(\frak{c}_{i, k+1})$,

\item[(c)] for $k < m$, $\frak{c}_{i, k} = \frak{c}_k$ does not depend on $i \in X$.

\item[(d)] for every $j < j'$ in $X$, $\dom(\frak{c}_{j, n-1}) < \dom(\frak{c}_{j', m})$ and

\item[(e)] $\otp(\dom(\frak{c}_{i, k}))$ does not depend on $i \in X$.

\end{itemize}

\end{claim}

Proof of Claim \ref{deltacr}: Just follow the argument in the proof of Claim \ref{CRccc} noting that $\dom(p_i)$'s are unbounded in $\omega_1$. \qed

\section{Countable joins}

In the course of club preservation arguments, we would like to be able to form new creatures out of old ones in the following way. Suppose $\langle q_i : i \geq 1 \rangle$ is a sequence of conditions in $\bb{Q} = \bb{Q}_{\CR, \Sigma}$ which forms a $\Delta$-system of an appropriate kind - It satisfies clauses (b)-(e) in Claim \ref{deltacr}. We'd like to construct a new condition $q \in \bb{Q}$ such that $q \Vdash_{\bb{Q}} ``\lim_n |\{i < n: q_i \in G_{\bb{Q}}\}|/n = 1$ and $\{i < \omega: q_i \notin G_{\bb{Q}}\}$ is infinite". This will require us to add ``countable joins" of certain sequences of creatures to $\CR$. This section introduces the countable join construction.

\begin{defn}

For $\alpha < \omega_1$, we say that $(\CR_p, \Sigma_p)$ is a partial $\aleph_1$-CP at $\alpha$ if for some $\aleph_1$-CP $(\CR, \Sigma)$,

\begin{itemize}

\item[(1)] $\CR_p = \CR \res \alpha = \{\frak{c} \in \CR: \sup(\dom(\frak{c})) < \alpha \}$ and

\item[(2)] $\Sigma_p = \Sigma \res \CR_p$.

\end{itemize} 

\end{defn}

\begin{defn}
\label{cp}
Suppose $k_{\star} \geq 1$, $\delta \in S_{k_{\star}}$ and $(\CR_p, \Sigma_p)$ is a partial $\aleph_1$-CP at $\delta$. Suppose $m < n < \omega$, $\delta \in S_{k_{\star}}$ for $k_{\star} \geq 1$ and $\bar{\frak{d}}_i = \langle \frak{d}_{i, k} : k < n \rangle$ satisfy the following for $1 \leq i < \omega$.

\begin{itemize}

\item[(a)] $\frak{d}_{i, k} \in \CR_p$.

\item[(b)] $\frak{d}_{i, j} = \frak{d}_j$ does not depend on $i$ for $j < m$.

\item[(c)] $\dom(\frak{d}_{i, k}) < \dom(\frak{d}_{i, k+1})$.

\item[(d)] $\dom(\frak{d}_{i, n - 1}) < \dom(\frak{d}_{i+1, m})$.

\item[(e)] $\otp(\dom(\frak{d}_{i, k}))$ only depends on $k$.

\item[(f)] $W = \bigcup \{\dom(\frak{d}_{i, k}) : 1 \leq i < \omega, k < n\}$ is unbounded in $\delta$ and has order type $\e + \omega^{j_{\star}}$ for some $\e < \omega_1$ and $1 \leq j_{\star} \leq k_{\star}$.

\end{itemize}

We say that $\langle \bar{\frak{d}}_i : i \geq 1 \rangle$ is a \textsf{joinable candidate} for $(\CR_p, \Sigma_p)$ at $\delta$. \\

For each $N \geq 1$ where $N$ is a power of $2$, we define new creatures $ \frak{c}^{\star}_{N} = (\dom(\frak{c}^{\star}_{N}), \pos(\frak{c}^{\star}_{N}), f_{\frak{c}^{\star}_{N}})$ and $\Sigma_{\star}(\frak{c}^{\star}_{N})$, as follows.

\begin{itemize}

\item[(1)] $\dom(\frak{c}^{\star}_1) = W$ and $\dom(\frak{c}^{\star}_N) = \bigcup \{\dom(\frak{d}_{i, k}) : N \leq i < \omega, m \leq k < n\}$ for $N \geq 2$.

\item[(2)] $f_{\frak{c}^{\star}_1} = \bigcup \{f_{\frak{d}_{i, k}}: 1 \leq i < \omega, k < n\}$ and $f_{\frak{c}^{\star}_N} = \bigcup \{f_{\frak{d}_{i, k}}: N \leq i < \omega, m \leq k < n\}$ for $N \geq 2$.

\item[(3)] $\Sigma_{\star}(\frak{c^{\star}_1})$ is the smallest family satisfying the following. 

\subitem (i) $\langle \frak{c}^{\star}_1 \rangle \in \Sigma_{\star}(\frak{c}^{\star}_1)$.

\subitem (ii) Whenever $j > 1$ is a power of $2$ and $\langle \frak{d}'_{i, k}: i < j, m \leq k < n \rangle$, $\langle \bar{\frak{f}}_{i, k}: i < j, m \leq k < n \rangle$ and $\langle \bar{\frak{g}}_k : k < m \rangle$ satisfy (a)-(d) below, we have, under appropriate order

$$\bigcup \{\bar{\frak{g}}_{k}: k < m\} \cup \bigcup \{\bar{\frak{f}}_{i, k} : i < j, m \leq k < n\} \cup \{\frak{c}^{\star}_{j}\} \in \Sigma_{\star}(\frak{c}^{\star}_1)$$

\subsubitem (a) $\frak{d}'_{i, k} \in \CR_p$ and $\dom(\frak{d}'_{i, k}) = \dom(\frak{d}_{i, k})$.

\subsubitem (b) $|\{i \in [j_1, j_2): (\exists k \in [m, n))(\frak{d}'_{i, k} \neq \frak{d}_{i, k}) \}| \leq (j_2 - j_1)/\log_2(j_1)$ for every $2 \leq j_1 < j_2 \leq j$ where $j_1, j_2$ are powers of $2$.

\subsubitem (c) $\bar{\frak{f}}_{i, k} \in \Sigma(\frak{d}'_{i, k})$.

\subsubitem (d) $\bar{\frak{g}}_k  \in \Sigma(\frak{d}_k)$.

\item[(4)] For $N \geq 2$, $\Sigma_{\star}(\frak{c}^{\star}_N)$ is the smallest family satisfying the following. 

\subitem (i) $\langle \frak{c}^{\star}_N \rangle \in \Sigma_{\star}(\frak{c}^{\star}_N)$.

\subitem (ii) Whenever $j > N$ is a power of $2$ and $\langle \frak{d}'_{i, k}: N \leq i < j, m \leq k < n \rangle$ and $\langle \bar{\frak{f}}_{i, k}: N \leq i < j, m \leq k < n \rangle$ satisfy (a)-(c) below, we have, under appropriate order

$$\bigcup \{\bar{\frak{f}}_{i, k} : i < j, m \leq k < n\} \cup \{\frak{c}^{\star}_{j}\} \in \Sigma_{\star}(\frak{c}^{\star}_N)$$

\subsubitem (a) $\frak{d}'_{i, k} \in \CR_p$ and $\dom(\frak{d}'_{i, k}) = \dom(\frak{d}_{i, k})$.

\subsubitem (b) $|\{i \in [j_1, j_2): (\exists k \in [m, n))(\frak{d}'_{i, k} \neq \frak{d}_{i, k}) \}| \leq (j_2 - j_1)/\log_2(j_1)$ for every $N \leq j_1 < j_2 \leq j$ where $j_1, j_2$ are powers of $2$.

\subsubitem (c) $\bar{\frak{f}}_{i, k} \in \Sigma(\frak{d}'_{i, k})$.

\item[(5)] $\pos(\frak{c}^{\star}_N) = \{\bigcup_{k < K} f_{\frak{c}_k}: \langle \frak{c}_{k}: k < K\rangle \in \Sigma_{\star}(\frak{c}^{\star}_N)\}$.

\end{itemize}

Let $(\CR'_p, \Sigma'_p)$ be the partial $\aleph_1$-CP at $\delta+1$ such that $\CR'_p = \CR_p \bigcup \{\frak{c}^{\star}_{N}:  N \geq 1 \text{ is a power of } 2 \}$ with $\dom(\frak{c}^{\star}_{N})$, $\pos(\frak{c}^{\star}_{N})$ and $f_{\frak{c}^{\star}_{N}}$ as above, $\Sigma'_p \res \CR_p = \Sigma_p$ and $\Sigma'_p(\frak{c}^{\star}_{N}) = \Sigma_{\star}(\frak{c}^{\star}_{N})$. We say that $(\CR'_p, \Sigma'_p)$ is the result of adding the countable join $\frak{c}_1 = \oplus_{i \geq 1} \bar{\frak{d}}_i$ of $\langle \bar{\frak{d}}_i : i \geq 1 \rangle$ to $(\Sigma_p, \CR_p)$.
\end{defn}

Note that $(\CR'_p, \Sigma'_p)$ is indeed a partial $\aleph_1$-CP at $\delta+1$ because $\Sigma'_p$ satisfies transitivity, cuts and finite joins.

\begin{lem}
\label{limit}
Let $(\CR'_p, \Sigma'_p)$ be as in Definition \ref{cp}. Let $(\CR, \Sigma)$ be an $\aleph_1$-CP such that $\CR'_p = \{\frak{c} \in \CR: \dom(\frak{c}) \subseteq \delta\}$ and $\Sigma'_p = \Sigma \res \CR'_p$. Let $\bb{Q} = \bb{Q}_{\CR, \Sigma}$, $p = \{\frak{c}^{\star}_1 = \oplus_{i \geq 1} \bar{\frak{d}}_i \}$ and $p_i = \{\frak{d}_{i, k} : k < n\}$. Then 
$$p \Vdash_{\bb{Q}} \lim_{j} \frac{|\{i < j: p_i \in G_{\bb{Q}}\}|}{j} = 1$$

\end{lem}

Proof of Lemma \ref{limit}: It suffices to show that for every $p_1 \geq p$ and $j_{\star} \geq 2^{10}$ there exists $p_2 \geq p_1$  such that 

$$p_2 \Vdash_{\bb{Q}} \frac{|\{i < j_{\star}: p_i \in \bb{G}_{\bb{Q}}\}|}{j_{\star}} > 1 - \frac{8}{\log_2{j_{\star}}}$$

Since $p_1 \geq p = \{\frak{c}^{\star}_1\}$, we can find $p_2 \geq p_1$ and $j_0 > j_{\star}$, such that $j_0$ is a power of $2$ and

$$\bigcup \{\bar{\frak{g}}_{k}: k < m\} \cup \bigcup \{\bar{\frak{f}}_{i, k} : i < j_0, m \leq k < n\} \cup \{\frak{c}^{\star}_{j_0}\} \subseteq p_2$$

where $\langle \frak{d}'_{i, k}: i < j_0, m \leq k < n \rangle$, $\langle \bar{\frak{f}}_{i, k}: i < j_0, m \leq k < n \rangle$ and $\langle \bar{\frak{g}}_k : k < m \rangle$ are as in Definition \ref{cp}(3)(ii). \\

Choose $N \geq 10$ such that $2^N \leq j_{\star} < 2^{N+1}$. Then $p_2$ forces that 

$$\frac{|\{i < j_{\star} : p_i \in G_{\bb{Q}}\}|}{j_{\star}} \geq 1 - \left( \sum_{1 \leq j < N} \frac{2^{j+1} - 2^j}{j j_{\star}} \right) - \frac{2^{N+1} - 2^N}{N j_{\star}} \geq 1 - \left( \sum_{1 \leq j < N}\frac{1}{j 2^{N - j}} \right) - \frac{1}{N} $$

Since $\sum_{1 \leq j < N/2} 1/(j 2^{N - j}) \leq N/2^{N/2} \leq 4/N$ (as $N \geq 10$) and $\sum_{N/2 \leq j < N} 1/(j 2^{N - j}) \leq 2/N$, it follows that 

$$p_2 \Vdash_{\bb{Q}} \frac{|\{i < j_{\star} : p_i \in G_{\bb{Q}}\}|}{j_{\star}} \geq 1 - \left(\frac{4}{N} + \frac{2}{N} + \frac{1}{N}\right) > 1 - \frac{8}{N}$$

 \qed

\begin{defn} 
\label{c2}
$(\CR, \Sigma)$ is a thin $\aleph_1$-CP if $(\CR, \Sigma)$ is an $\aleph_1$-CP and there exist $S$ and $\langle \frak{c}_{\delta}: \delta \in S\rangle$ such that the following hold.

\begin{itemize}

\item[(a)] $S \subseteq \bigcup_{k \geq 1} S_k$.

\item[(b)] $\frak{c}_{\delta} \in \CR$.

\item[(c)] For every $k_{\star} \geq 1$ and $\delta \in S \cap S_{k_{\star}}$, letting $(\CR_P, \Sigma_p)$ be the partial $\aleph_1$-CP at $\delta$ satisfying $\CR_p = \CR \res \delta = \{\frak{c} \in \CR: \sup(\dom(\frak{c})) < \delta\}$ and $\Sigma_p = \Sigma \res \CR_p$, there exists a joinable candidate $\langle \bar{\frak{d}}_i : i \geq 1 \rangle$ for $(\CR_p, \Sigma_p)$ at $\delta$ such that

\subitem (i) $\frak{c}_{\delta} = \oplus_{i \geq 1} \bar{\frak{d}}_i$ and

\subitem (ii) $\CR'_p = \{\frak{c} \in \CR: \dom(\frak{c}) \subseteq \delta\}$ and $\Sigma'_p = \Sigma \res \CR'_p$ where $(\CR'_p, \Sigma'_p)$ is the result of adding $\oplus_{i \geq 1} \bar{\frak{d}}_i$ to $(\CR_p, \Sigma_p)$.

\item[(d)] $\frak{c} \in \CR$ iff $\frak{c}$ is a finite join of $\{\frak{d} \in \CR: \frak{d} \text{ is finite} \} \cup \bigcup \{\Sigma(\frak{c}_{\delta}): \delta \in S\}$.
 
\end{itemize}

\end{defn}

\begin{claim}
\label{tailagree}
Suppose $(\CR, \Sigma)$ is an $\aleph_1$-CP as witnessed by $S, \langle \frak{c}_{\delta} : \delta \in S\rangle$. Suppose $\frak{c} \in \CR$, $k_{\star} \geq 1$, $\delta \in S_{k_{\star}}$, $\dom(\frak{c})$ is an unbounded subset of $\delta$. Then there exist $\bar{c} = \langle \frak{c}_k : k \leq k_1 \rangle \in \Sigma(\frak{c})$ and $\bar{d} = \langle \frak{d}_k : k \leq k_2 \rangle \in \Sigma(\frak{c}_{\delta})$ such that $\frak{c}_{k_1} = \frak{d}_{k_2}$.

\end{claim}

Proof of Claim \ref{tailagree}: Easily follows from Definition \ref{c2}. \qed

\section{Guided products}

\begin{defn}
\label{gprod}
Suppose $\langle \bb{Q}_{\alpha}: \alpha < \omega_2 \rangle$ and $\langle p^{\star}_{\delta}: \delta \in S^{\omega_2}_{\aleph_0} \rangle$ satisfy the following.

\begin{itemize}

\item[(i)] $\bb{Q}_{\alpha} = \bb{Q}_{\CR_{\alpha}, \Sigma_{\alpha}}$ where $(\CR_{\alpha}, \Sigma_{\alpha})$ is a thin $\aleph_1$-CP.

\item[(ii)] $p^{\star}_{\delta}$ is a function whose domain is a countable  unbounded subset of $\delta$ and for every $\alpha \in \dom(p^{\star}_{\delta})$, $p^{\star}_{\delta}(\alpha) \in \bb{Q}_{\alpha}$.

\end{itemize}

For $\gamma \leq \omega_2$, define a forcing $\bb{P}_{\gamma}$ as follows.

\begin{itemize}

\item[(1)] $p \in \bb{P}_{\gamma}$ iff 

\subitem (a) $p$ is a function, $\dom(p)$ is a countable subset of $\gamma$,

\subitem (b) for every $\alpha \in \dom(p)$, $p(\alpha) \in \bb{Q}_{\alpha}$ and

\subitem (c) for every $\delta \leq \gamma$ with $\cf(\delta) = \aleph_0$, if $\dom(p) \cap \delta$ is unbounded in $\delta$, then for some $\eta < \delta$, $p \res (\eta, \delta) = p^{\star}_{\delta} \res (\eta, \delta)$.

\item[(2)] For $p, q \in \bb{P}_{\gamma}$, define $p \leq q$ iff $\dom(p) \subseteq \dom(q)$ and for every $\alpha \in \dom(p)$, $p(\alpha) \leq_{\bb{Q}_{\alpha}} q(\alpha)$.

\end{itemize}

We say that $\bb{P}_{\omega_2}$ is the countable support product of $\langle \bb{Q}_{\alpha}: \alpha < \omega_2 \rangle$ guided by $\langle p^{\star}_{\delta}: \delta \in S^{\omega_2}_{\aleph_0} \rangle$. Note that for $\cf(\gamma) = \aleph_1$, $\bb{P}_{\gamma}$ is completely determined by $\langle \bb{Q}_{\alpha}: \alpha < \gamma \rangle$ and $\langle p^{\star}_{\delta}: \delta < \gamma, \cf(\delta) = \aleph_0 \rangle$.

\end{defn}

\begin{claim} 
\label{cc2}
Let $\langle \bb{Q}_{\alpha}: \alpha < \omega_2 \rangle$, $\langle p^{\star}_{\delta}: \delta \in S^{\omega_2}_{\aleph_0} \rangle$ and $\bb{P}_{\gamma}$ for $\gamma \leq \omega_2$ be as in Definition \ref{gprod}. Then the following hold.

\begin{itemize}

\item[(a)] $\bb{P}_{\gamma + 1} = \bb{P}_{\gamma} \times \bb{Q}_{\gamma}$.

\item[(b)] $\bb{P}_{\gamma}$ satisfies ccc.

\end{itemize}

\end{claim}

Proof of Claim \ref{cc2}: (a) is obvious from the definition of $\bb{P}_{\gamma}$. (b) follows from Claim \ref{deltait} below and the fact that each $\bb{Q}_{\alpha}$ has $\aleph_1$ as a precaliber (Claim \ref{CRccc}).

\begin{claim}
\label{deltait}
Suppose $\gamma \leq \omega_2$ and $\langle p_i : i < \omega_1 \rangle$ is a sequence of conditions in $\bb{P}_{\gamma}$. Then there exists $X \in [\omega_1]^{\aleph_1}$ and a finite $F \subseteq \omega_2$ such that for every $\alpha \in \omega_2 \setminus F$, if there are $i < j$ in $X$ such that $\alpha \in \dom(p_i) \cap \dom(p_j)$, then $(\forall i \in X)(\alpha \in \dom(p_i)$ and $p_i(\alpha)$ does not depend on $i \in X)$.

\end{claim}

Proof of Claim \ref{deltait}: By induction on $\gamma \leq \omega_2$. If $\gamma$ is a successor or $\gamma = \omega_2$, this is trivial. \\

Suppose $\cf(\gamma) = \aleph_0$ and let $\langle p_i : i < \omega_2 \rangle$ be a sequence of conditions in $\bb{P}_{\gamma}$. Let $\langle \gamma_n : n < \omega \rangle$ be increasing cofinal in $\gamma$. For each $i < \omega_1$, choose $n = n_i < \omega$ such that either $p_i \in \bb{P}_{\gamma_n}$ or $p_i \res (\gamma_n, \gamma) = p^{\star}_{\gamma} \res (\gamma_n, \gamma)$. Choose $X \in [\omega_1]^{\aleph_1}$ and $n_{\star} < \omega$, such that $(\forall i \in X)(n_1 = n_{\star})$ and apply the inductive hypothesis to $\langle p_i \res \gamma_{n_{\star}}: i \in X \rangle$. \\

Next suppose $\cf(\gamma) = \omega_1$ and let $\langle p_i : i < \omega_2 \rangle$ be a sequence of conditions in $\bb{P}_{\gamma}$. Choose $\langle \gamma_i : i < \omega_1 \rangle$ continuously increasing and cofinal in $\gamma$ such that $\cf(\gamma_i) = \aleph_0$ for every $i < \omega_1$. For each $i < \omega_1$, choose $j = j_i < i$ such that either $p_i \res \gamma_i \in \bb{P}_{\gamma_j}$ or $p_i \res (\gamma_j, \gamma_i) = p^{\star}_{\gamma_i} \res (\gamma_j, \gamma_i)$. By Fodor's lemma we can get $S \in [\omega_1]^{\aleph_1}$ and $j_{\star} < \omega_1$ such that $(\forall i \in S)(j_i = j_{\star})$. Choose $X \in [S]^{\aleph_1}$ such that for every $i < j$ in $X$, $\dom(p_i) \cap \dom(p_j) \subseteq \gamma_{j_{\star}}$. Now apply the inductive hypothesis to $\langle p_i \res \gamma_{j_{\star}}: i \in S \rangle$. \qed

\begin{lem}
\label{clubstar}
Let $\langle \bb{Q}_{\alpha}: \alpha < \omega_2 \rangle$, $\langle p^{\star}_{\delta}: \delta \in S^{\omega_2}_{\aleph_0} \rangle$ and $\bb{P}_{\omega_2}$ be as in Definition \ref{gprod}. Then $V^{\bb{P}_{\omega_2}} \models \neg \club^1$. 
\end{lem}

Proof of Lemma \ref{clubstar}: Towards a contradiction, suppose $p_0 \in \bb{P}_{\omega_2}$, $\langle \name{A}_{\delta} = \{\name{\alpha}_{\delta, n} : n < \omega\} : \delta \in \Lim(\omega_1) \rangle \in V^{\bb{P}_{\omega_2}}$ are such that $p_0 \Vdash ``(\forall \delta \in \Lim(\omega_1))(\{\name{\alpha}_{\delta, n} : n < \omega\}$ is increasing cofinal in $\delta)$ and $\langle \name{A}_{\delta} : \delta < \omega_1 \rangle$ is a $\club^1$ witnessing sequence". Since $\bb{P}_{\omega_2}$ satisfies ccc, we can find $\gamma < \omega_2$ such that $p_0 \in \bb{P}_{\gamma}$ and each $\name{\alpha}_{\delta, n}$ is a $\bb{P}_{\gamma}$-name. \\

Let $\name{X} = \{\alpha < \omega_1: \name{f}_{\bb{Q}_{\gamma}} = 1\}$. Then $\name{X} \in V^{\bb{P}_{\gamma + 1}}$ and $ V^{\bb{P}_{\gamma + 1}} \models \name{X} \in [\omega_1]^{\aleph_1}$. So there exist $p_1 \in \bb{P}_{\gamma}$, $q \in \bb{Q}_{\gamma}$, $\delta \in \Lim(\omega_1)$ and $n_{\star} < \omega$ such that $p_1 \geq p_0$ and $(p_1, q) \Vdash_{\bb{P}_{\gamma + 1}} (\forall n \geq n_{\star})(\name{\alpha}_{\delta, n} \in \name{X})$. Note that we must have that $\dom(q) \cap \delta$ is unbounded in $\delta$ otherwise we can easily extend $(p_1, q)$ to get a contradiction. By possibly extending $q$, by Definition \ref{crdefn}(D)(ii), we can assume that $q = \{\frak{c}_k : k < K_{\star}\}$ where $\dom(\frak{c}_k) < \dom(\frak{c}_{k+1})$ for every $k < K_{\star}-1$ and for some $K < K_{\star}$, $\dom(\frak{c}_{K})$ is an unbounded subset of $\delta$. Let $S_{\gamma}$ and $\langle \frak{c}_{\gamma, \delta}: \delta \in S_{\gamma} \rangle$ witness that $(\CR_{\gamma}, \Sigma_{\gamma})$ is a thin $\aleph_1$-CP. By Claim \ref{tailagree}, we can further assume that $\frak{c}_K = \frak{c}'_{K'}$ for some $\langle \frak{c}'_n : n \leq K' \rangle \in \Sigma(\frak{c}_{\gamma, \delta})$. \\

Let $m < n < \omega$ and $\bar{\frak{d}}_i = \langle \frak{d}_{i, k} : k < n \rangle$ for $i \geq 1$ be as in Definition \ref{cp} and $\frak{c}_{\gamma, \delta} = \oplus_{i \geq 1} \bar{\frak{d}}_i$. Then as $\langle \frak{c}'_n : n \leq K' \rangle \in \Sigma(\frak{c}_{\gamma, \delta})$, we can find $N \geq 1$ a power of $2$ such that $\frak{c}_K = \frak{c}'_{K'} = \frak{c}^{\star}_N$ in the notation of Definition \ref{cp}. \\

Choose $p_2 \in \bb{P}_{\gamma}$, $p_2 \geq p_1$, $n(1) > n_{\star}$ and $\alpha > \min(\dom(\frak{c}^{\star}_N))$ such that $p_2 \Vdash_{\bb{P}_{\gamma}} \name{\alpha}_{\delta, n(1)} = \alpha$. We can assume that $\alpha \in \dom(\frak{c}^{\star}_N)$ - Otherwise letting $\frak{c}_{\star}$ be a creature with domain $\{\alpha\}$ and $f_{\frak{c}_{\star}}(\alpha) = 1$, we have $q' = q \cup \{\frak{c}_{\star}\} \Vdash_{\bb{Q}_{\gamma}} \alpha \notin \name{X}$ so that $(p_2, q')$ forces a contradiction. Choose $i(\star) \geq N$ and $m \leq k(\star) < n$ such that $\alpha \in \dom(\frak{d}_{i(\star), k})$. Let $N_1 \geq N$ be the largest power of $2$ such that $i(\star) \geq N_1$ and let $j > i(\star)$ be a power of $2$. Choose a creature $\frak{d}'_{i(\star), k(\star)}$ such that $\dom(\frak{d}'_{i(\star), k(\star)}) = \dom(\frak{d}_{i(\star), k(\star)})$ and $\bar{\frak{f}} \in \Sigma(\frak{d}'_{i(\star), k(\star)})$ such that for some finite $\frak{c}_{\star} \in \bar{\frak{f}}$, $\dom(\frak{c}_{\star}) = \{\alpha\}$ and $f_{\frak{c}_{\star}}(\alpha) = 0$. It follows that, under appropriate order

$$\{\frak{d}_{i, k}: N \leq i < j, m \leq k < n, (i, k) \neq (i(\star), k(\star))\} \cup \bar{\frak{f}} \cup \{\frak{c}^{\star}_j\} \in \Sigma(\frak{c}^{\star}_N)$$

Let $q' = (q \setminus \{\frak{c}_K\}) \cup \{\frak{d}_{i, k}: N \leq i < j, m \leq k < n, (i, k) \neq (i(\star), k(\star))\} \cup \bar{\frak{f}} \cup \{\frak{c}^{\star}_j\}$. Then $(p_2, q') \geq (p, q)$ and $q' \Vdash_{\bb{Q}_{\gamma}} \alpha \notin \name{X}$ - Contradiction. \qed

\section{$\club^{\lim}$ and $\neg \club^1$}

We define a preparatory forcing $\bb{R}$ which generically adds $\langle \bb{Q}_{\alpha}: \alpha < \omega_2 \rangle$ and $\langle p^{\star}_{\delta}: \delta \in S^{\omega_2}_{\aleph_0} \rangle$ satisfying Definition \ref{gprod}(i)-(ii) using countable approximations. This ensures that the resulting guided product $\bb{P}_{\omega_2}$ preserves a $\club^{\lim}$ witnessing sequence $\bar{A}$ which is also added by $\bb{R}$ via countable approximations.

\begin{defn}
Let $\bb{R}$ be a forcing whose conditions are $r = (u_r, \delta_r, \langle \bb{Q}_{r, \alpha}: \alpha \in u_r \rangle, v_r, \langle p^{\star}_{r, \alpha} : \alpha \in v_r \rangle, \bar{A}_{r})$ where

\begin{itemize}

\item[(a)] $u_r \in [\omega_2]^{\leq \aleph_0}$, $\delta_r < \omega_1$,

\item[(b)] $\bb{Q}_{r, \alpha} = \bigcup_{\xi < \delta_r} (\bb{Q}_{\CR_{r, \alpha}, \Sigma_{r, \alpha}} \res \xi)$ for some thin $\aleph_1$-CP $(\CR_{r, \alpha}, \Sigma_{r, \alpha})$ as witnessed by $(S_{r, \alpha}, \langle \frak{c}_{r, \alpha, \delta}: \delta \in S_{r, \alpha} \rangle)$ - So only $S_{r, \alpha} \cap \delta_r$ and $\langle \frak{c}_{r, \alpha, \delta}: \delta \in S_{r, \alpha} \cap \delta_r \rangle$ are relevant,

\item[(c)] $v_r \subseteq u \cap \Lim(\omega_1)$ and for every $\alpha \in v_r$, $u_r \cap \alpha$ is unbounded in $\alpha$,

\item[(d)] $p^{\star}_{r, \alpha}$ is a function with domain an unbounded subset of $u_r \cap \alpha$ and for each $\xi \in \dom(p^{\star}_{r, \alpha})$, $p^{\star}_{r, \alpha}(\xi) \in \bb{Q}_{r, \alpha}$ and

\item[(e)] $\bar{A}_r = \langle A_{r, \gamma}: \gamma \in \Lim(\omega_1) \cap \delta_r \rangle$ where each $A_{r, \gamma}$ is an unbounded subset of $\gamma$ of order type $\omega$.

\end{itemize}

For $r, s \in \bb{R}$, define $r \leq s$ iff the following hold.

\begin{itemize}

\item[(i)] $u_r \subseteq u_s$, $\delta_r \leq \delta_s$.

\item[(ii)] For every $\alpha \in u_r$, $S_{r, \alpha} \cap \delta_r = S_{s, \alpha} \cap \delta_r$ and $\frak{c}_{r, \alpha, \delta} = \frak{c}_{s, \alpha, \delta}$ for every $\delta \in S_r \cap \delta_r$. It follows that $\bb{Q}_{r, \alpha} \subseteq \bb{Q}_{s, \alpha}$ and for every $p \in \bb{Q}_{s, \alpha}$, if $\dom(p)$ is bounded below $\delta_r$, then $p \in \bb{Q}_{r, \alpha}$.

\item[(iii)] $v_r \subseteq v_s$ and for every $\alpha \in v_r$, $p^{\star}_{s, \alpha} = p^{\star}_{r, \alpha}$.

\item[(iv)] $\bar{A}_r = \bar{A}_s \res (\Lim(\omega_1) \cap \delta_r)$.

\end{itemize}

\end{defn}

\begin{claim}
\label{onR}
$\bb{R}$ is countably closed and hence it preserves stationary subsets of $\omega_1$. Under CH, it satisfies $\aleph_2$-c.c. and therefore preserves all cofinalities.

\end{claim}

Proof of Claim \ref{onR}: It is clear that $\bb{R}$ is countably closed. Next let $\{r_i : i < \omega_2\} \subseteq \bb{R}$. Using CH, we can find $X_0 \in [\omega_2]^{\aleph_2}$ such that $\langle u_{r_i}: i \in X_0 \rangle$ forms a $\Delta$-system with root $u_{\star}$. By possibly extending each $r_i$, we can assume that $u_{r_i} \setminus u_{\star} \neq \emptyset$ for every $i \in X_0$. Choose $X \in [X_0]^{\aleph_2}$ such that the following hold.

\begin{itemize}

\item[(i)] For every $i, j \in X$ with $i < j$, $\sup(u_{\star}) < \min(u_{r_i} \setminus u_{\star}) \leq \sup(u_{r_i} \setminus u_{\star}) < \inf(u_{r_j} \setminus u_{\star})$.

\item[(ii)] $\langle v_{r_i}: i \in X \rangle$ forms a $\Delta$-system with root $v_{\star} \subseteq u_{\star}$ .

\item[(iii)] $\delta_{r_i} = \delta_{\star}$ does not depend on $i \in X$.

\item[(iv)] For every $\alpha \in u_{\star}$, $\bb{Q}_{r_i, \alpha} = \bb{Q}_{\alpha}$ does not depend on $i \in X$.

\item[(v)] For every $\alpha \in v_{\star}$, $p^{\star}_{r_i, \alpha} = p^{\star}_{\alpha}$ does not depend on $i \in X$.

\item[(vi)] $\bar{A}_{r_i} = \bar{A}_{\star}$ does not depend on $i \in X$.

\end{itemize}

For clauses (iv), (v) and (vi), we use CH. It is clear that any two conditions in $\{r_i : i \in X\}$ have a common extension. \qed \\

\relax From now on we assume CH. The next claim is easily verified.

\begin{claim} Each of the following sets is dense in $\bb{R}$.

\begin{itemize}

\item[(a)] $\{r \in \bb{R}: \alpha \in u_r\}$ for $\alpha < \omega_2$.

\item[(b)] $\{r \in \bb{R}: \delta_r > \delta\}$ for $\delta < \omega_1$.

\item[(c)] $\{r \in \bb{R}: \delta \in v_r\}$ for $\delta \in S^{\omega_2}_{\aleph_0}$.

\end{itemize}

\end{claim}

Let $G_{\bb{R}}$ be $\bb{R}$-generic over $V$. Work in $V_1 = V[G_{\bb{R}}]$. For $\alpha < \omega_2$, define $\bb{Q}_{\alpha} = \bigcup \{\bb{Q}_{r, \alpha}: r \in G_{\bb{R}}, \alpha \in u_r\}$. Note that for every $\alpha < \omega_2$, $S_{\alpha} = \bigcup \{S_{r, \alpha} \cap \delta_r: r \in G_{\bb{R}}, \alpha \in u_r\}$ is a stationary subset of $\bigcup_{k \geq 1} S_k$ and $V_1 \models ``(\forall \alpha < \omega_2)(\bb{Q}_{\alpha} = \bb{Q}_{\CR_{\alpha}, \Sigma_{\alpha}}$ for some thin $\aleph_1$-CP $(\CR_{\alpha}, \Sigma_{\alpha}))$". For $\delta \in S^{\omega_2}_{\aleph_0}$, let $p^{\star}_{\delta} = p^{\star}_{r, \delta}$ for some $r \in G_{\bb{R}}$ with $\delta \in v_r$. Let $\bar{A} = \langle A_{\delta}: \delta \in \Lim(\omega_1) \rangle = \bigcup \{\bar{A}_r : r \in G_{\bb{R}}\}$. Let $\{\alpha_{\delta, n}: n < \omega\}$ list $A_{\delta}$ in increasing order.  \\

Let $\bb{P}_{\omega_2} \in V_1$ be the countable support product of $\langle \bb{Q}_{\alpha}: \alpha < \omega_2 \rangle$ guided by $\langle p^{\star}_{\delta} : \delta \in S^{\omega_2}_{\aleph_0} \rangle$. Note that, since $\bb{R}$ is countably closed, the set of conditions $(r, p) \in \bb{R} \star \bb{P}_{\omega_2}$ satisfying the following is dense in $\bb{R} \star \bb{P}_{\omega_2}$.

\begin{itemize}

\item[(a)] $p$ is an actual object.

\item[(b)] $\dom(p) \subseteq u_r$.

\item[(c)] $(\forall \alpha \in \dom(p))(p(\alpha) \in \bb{Q}_{r, \alpha})$.

\item[(d)] For every $\alpha <  \omega_2$ of cofinality $\aleph_0$, if $\dom(p) \cap \alpha$ is unbounded in $\alpha$, then $\alpha \in v_r$.

\end{itemize}

So we can assume that our conditions in $\bb{R} \star \bb{P}_{\omega_2}$ have this form.

\begin{thm} 
\label{main}
$V_1^{\bb{P}_{\omega_2}} \models \club^{\lim} \wedge \neg \club^1$

\end{thm}

Proof of Theorem \ref{main}: That $V_1^{\bb{P}_{\omega_2}} \models \neg \club^1$ follows from Lemma \ref{clubstar}. We'll show that $\bar{A}$ witnesses $\club^{\lim}$ in $V_1^{\bb{P}_{\omega_2}}$. Suppose $(r_{\star}, p_{\star}) \Vdash_{\bb{R} \star \bb{P}_{\omega_2}} \name{A} \in [\omega_1]^{\aleph_1}$. We'll construct $(r, p) \geq (r_{\star}, p_{\star})$ and $\delta < \omega_1$ such that 

$$(r, p) \Vdash_{\bb{R} \star \bb{P}_{\omega_2}} \lim_n \frac{|\{k < n: \name{\alpha}_{\delta, k} \in \name{A}\}|}{n} = 1$$ \\

Choose $\langle (r_i, p_i, \gamma_i): i < \omega_1 \rangle$ such that the following hold.

\begin{itemize}

\item[(i)] $(r_i, p_i) \geq (r_{\star}, p_{\star})$.

\item[(ii)] For all $i < j < \omega_1$, $r_i \leq_{\bb{R}} r_j$, $\sup(u_{r_i}) < \sup(u_{r_j})$ and $i \leq \delta_{r_i} < \delta_{r_j}$. 

\item[(iii)] For every $i < \omega_1$, $\sup(\bigcup_{j < i} \dom(p_j)) < \sup(\dom(p_i))$. 

\item[(iv)] For every $i < \omega_1$ $i \in u_i$ and for every $\alpha < \sup(u_{r_i})$, there exists $ j \in (i, \omega_1)$ such that $\alpha \in u_{r_j}$. So $\bigcup_{i < \omega_1} u_{r_i} =\alpha_{\star} \in [\omega_1, \omega_2)$ and $\cf(\alpha_{\star}) = \aleph_1$.

\item[(v)] For every $\delta <\alpha_{\star}$ with $\cf(\delta) = \aleph_0$, there exists $i < \omega_1$ such that $\delta \in v_{r_i}$. Hence $\bigcup_{i < \omega_1} v_{r_i} = \{\delta < \alpha_{\star}: \cf(\delta) = \aleph_0\}$.

\item[(vi)] $\langle \gamma_i : i < \omega_1 \rangle$ is a strictly increasing sequence in $\omega_1$.

\item[(vii)] $(r_i, p_i) \Vdash \gamma_i \in \name{A}$.

\end{itemize}

\begin{claim}
\label{deltait1}
There exist $F \subseteq \omega_2$ finite and $X \in [\omega_1]^{\aleph_1}$ such that for every $\alpha \in \omega_2 \setminus F$, if $\alpha \in \dom(p_i) \cap \dom(p_j)$ for some $i < j$ in $X$, then $(\forall i \in X)(\alpha \in \dom(p_i)$ and $p_i(\alpha)$ does not depend on $i \in X)$. 

\end{claim}

Proof of Claim \ref{deltait1}: For $\alpha < \alpha_{\star}$, let $\bb{Q}'_{\alpha} = \bigcup \{\bb{Q}_{r_i, \alpha}: i < \omega_1, \alpha \in u_{r_i}\}$. Then $\bb{Q}'_{\alpha}$ is a thin $\aleph_1$-CP. For $\delta < \alpha_{\star}$ with $\cf(\delta) = \aleph_0$, let $p^{\star}_{\delta} = p^{\star}_{r_i, \delta}$ where $i < \omega_1$ and $\delta \in v_{r_i}$. Let $\bb{P}_{\alpha_{\star}}$ be the countable support product of $\langle \bb{Q}'_{\alpha} : \alpha < \alpha_{\star} \rangle$ guided by $\langle p^{\star}_{\delta}: \delta < \alpha_{\star}, \cf(\delta) = \aleph_0\rangle$ so that each $p_i \in \bb{P}_{\alpha_{\star}}$. Now apply Claim \ref{deltait}. \qed \\

By shrinking $X$ and $F$, we can assume that for every $i \in X$, $F \subseteq \dom(p_i)$. Let $W = \bigcap_{i \in X} (\dom(p_i) \setminus F)$ and $Y_i = \dom(p_i) \setminus (F \cup W)$. Then $\langle Y_i: i \in X \rangle$ is a sequence of pairwise disjoint non empty countable sets. By shrinking $X$, we can also assume that for every $i < j$ in $X$, $\sup(Y_i) < \min(Y_j)$ and $\otp(\dom(p_i))$ does not depend on $i \in X$. \\

By Claim \ref{deltacr}, we can find $X_1 \in [X]^{\aleph_1}$ such that for every $\alpha \in F$ exactly one of the following holds.

\begin{itemize}

\item[(A)] For every $i \in X_1$, $p_i(\alpha) = q_{\alpha}$ does not depend on $i$.

\item[(B)] There are $m = m_{\alpha}$, $n = n_{\alpha}$, $m < n < \omega$ and $\langle q_{i, \alpha}: i \in X_1 \rangle$ such that for every $i \in X_1$,

\subitem (i) $q_{i, \alpha} \in \bb{Q}_{r_i, \alpha}$, $\dom(q_{i, \alpha}) = \dom(p_i(\alpha))$ and $r_i \Vdash_{\bb{R}} p_i(\alpha) \leq_{\bb{Q}_{\alpha}} q_{i, \alpha}$,

\subitem (ii) $q_{i, \alpha} = \{\frak{d}_{i, \alpha, k}: k < n\}$ and for every $k < n-1$, $\dom(\frak{d}_{i, \alpha, k}) < \dom(\frak{d}_{i, \alpha, k+1})$,

\subitem (iii) for every $k < m$, $\frak{d}_{i, \alpha, k} = \frak{d}_{\alpha, k}$ does not depend on $i \in X_1$,

\subitem (iv) for every $j < j'$ in $X$, $\dom(\frak{d}_{j, \alpha, n-1}) < \dom(\frak{d}_{j', \alpha, m})$ and

\subitem (v) $\otp(\frak{d}_{i, \alpha, k}) = \theta_{\alpha, k}$ does not depend on $i \in X_1$ and $1 \leq k_{\alpha} < \omega$ is such that $\theta_{\alpha, k} < \omega^{k_{\alpha}}$.

\end{itemize}

Let $F_0$ be the set of $\alpha \in F$ for which case (A) holds and $F_1 = F \setminus F_0$. \\

By reindexing, we can assume that $X_1 = \omega_1$. Let $k_{\star} = \max (\{k_{\alpha} + 2: \alpha \in F\})$. Put $Y = \bigcup_{i < \omega_1} Y_i$. Choose a club $E \subseteq \omega_1$ such that for every $\delta \in E$, the following hold.

\begin{itemize}

\item[(a)] For every $i < \delta$, there exists $j < \delta$ such that $\sup(u_{r_i} \cap Y) < \sup(Y_j)$.

\item[(b)] $\sup(\{\delta_{r_i}: i < \delta\}) = \delta$.

\item[(c)] For every $\alpha \in F_1$, $\sup(\{\dom(q_i(\alpha)): i < \delta\}) = \delta$.

\item[(d)] $\sup(\{\gamma_i : i < \delta\}) = \delta$.

\end{itemize}

Fix $\delta \in S_{k_{\star}} \cap E$ and let $\langle i(n): n < \omega \rangle$ be increasing cofinal in $\delta$. Let $\alpha_{\star} = \sup(\{Y_{i(n)}: n < \omega\})$. We can assume that $\alpha_{\star} \notin F \cup W$ - Just pick a sufficiently large $\delta \in S_{k_{\star}} \cap E$. Define $r \in \bb{R}$ as follows.

\begin{itemize}

\item[(a)] $u_r = \bigcup_{n < \omega} u_{r_{i(n)}} \cup \{\alpha_{\star}\}$, $\delta_r = \delta + 1$.

\item[(b)] For $\alpha \in u_r$, choose $\bb{Q}_{r, \alpha}$, $(\CR_{r, \alpha}, \Sigma_{r, \alpha})$ and $(S_{r, \alpha}, \langle \frak{c}_{r, \alpha, \delta}: \delta \in S_{r, \alpha} \rangle)$ as follows.

\subitem (i) If $\alpha \in u_r \setminus (F_1 \cup \{\alpha_{\star}\})$, choose a thin $\aleph_1$-CP $(\CR_{r, \alpha}, \Sigma_{r, \alpha})$ with witnessing pair $(S_{r, \alpha}, \langle \frak{c}_{r, \alpha, \delta}: \delta \in S_{r, \alpha} \rangle)$ such that for every $n < \omega$, $S_{r, \alpha} \cap \delta_{i(n)} = S_{r_{i(n)}, \alpha} \cap \delta_{i(n)}$ and $\frak{c}_{r, \alpha, \delta} = \frak{c}_{r_{i(n)}, \alpha, \delta}$ for every $\delta \in S_{r, \alpha} \cap \delta_{i(n)}$. So $\bigcup_{n < \omega} \bb{Q}_{r_{i(n)}, \alpha} \subseteq \bb{Q}_{r, \alpha} = \bb{Q}_{\CR_{r, \alpha}, \Sigma_{r, \alpha}} \res \delta$.

\subitem (ii) If $\alpha = \alpha_{\star}$, choose $\bb{Q}_{r, \alpha}$, $(\CR_{r, \alpha}, \Sigma_{r, \alpha})$ and $(S_{r, \alpha}, \langle \frak{c}_{r, \alpha, \delta}: \delta \in S_{r, \alpha} \rangle)$ arbitrarily.

\subitem (iii) If $\alpha \in F_1$, choose a thin $\aleph_1$-CP $(\CR_{r, \alpha}, \Sigma_{r, \alpha})$ with witnessing pair $(S_{r, \alpha}, \langle \frak{c}_{r, \alpha, \delta}: \delta \in S_{r, \alpha} \rangle)$ such that for every $n < \omega$, $S_{r, \alpha} \cap \delta_{i(n)} = S_{r_{i(n)}, \alpha} \cap \delta_{i(n)}$, $\frak{c}_{r, \alpha, \delta} = \frak{c}_{r_{i(n)}, \alpha, \delta}$ for every $\delta \in S_{r, \alpha} \cap \delta_{i(n)}$, $\delta \in S_{r, \alpha}$ and $\frak{c}_{r, \alpha, \delta} = \oplus_{n \geq 1} \langle \frak{d}_{i(n), \alpha, k}: k < n_{\alpha} \rangle$ where $\langle \frak{d}_{i(n), \alpha, k}: k < n_{\alpha} \rangle$ is from clause (B)(ii) above. Put $\bb{Q}_{r, \alpha} = \bb{Q}_{\CR_{r, \alpha}, \Sigma_{r, \alpha}} \res \delta$.

\item[(c)] $v_r = \bigcup_{n < \omega} v_{r_{i(n)}} \cup \{\alpha_{\star}\}$.

\item[(d)] For $\alpha \in v_{r_{i(n)}}$, $p^{\star}_{r, \alpha} = p^{\star}_{r_{i(n)}, \alpha}$ and $p^{\star}_{r, \alpha_{\star}} = \bigcup_{n < \omega} p_{i(n)} \res Y_{i(n)}$. So $\dom(p^{\star}_{r, \alpha_{\star}})$ is an unbounded subset of $u_r \cap \alpha_{\star}$.

\item[(e)] $\bar{A}_{r} = \bigcup_{n < \omega} \bar{A}_{r_{i(n)}} \cup \{(\delta, \{\gamma_{i(n)}: n < \omega\})\}$.

\end{itemize}

Next define $p$ as follows. 

\begin{itemize}

\item[(i)] $\dom(p) = F \cup W \cup \bigcup_{n < \omega} Y_n$.

\item[(ii)] If $\alpha \in F_0$, then $p(\alpha) = q_{\alpha}$ where $q_{\alpha}$ is from clause (B) above.

\item[(iii)] If $\alpha \in F_1$, then $p(\alpha) = \{\frak{c}_{r, \alpha, \delta}\}$.

\item[(iv)] If $\alpha \in W$, then $p(\alpha) = p_{i(n)}(\alpha)$ which does not depend on $n < \omega$.

\item[(v)] For every $n < \omega$, $p \res Y_{i(n)} = p_{i(n)} \res Y_i$.

\end{itemize}

It is clear that $(r_{\star}, p_{\star}) \leq_{\bb{R} \star \bb{P}_{\omega_2}} (r, p)$. By Lemma \ref{limit}, 

$$(r, p) \Vdash_{\bb{R} \star \bb{P}_{\omega_2}} \lim_n \frac{|\{k < n: (r_{i(k)}, p_{i(k)}) \in G_{\bb{R} \star \bb{P}_{\omega_2}}\}|}{n} = 1$$

Hence

$$(r, p) \Vdash_{\bb{R} \star \bb{P}_{\omega_2}} \lim_n \frac{|\{k < n: \gamma_{i(k)} \in \name{A}\}|}{n} = 1$$

Since $A_{r, \delta} = \{\gamma_{i(n)}: n < \omega\}$, the result follows.  \qed

\section{On $\club^{\inf \geq a}$}

\begin{defn}
\label{definfminus}
For $a \in (0, 1]$, the principle $\club^{\inf > a -}$ says the following. There exists $\bar{A} = \langle A_{\delta} : \delta \in \Lim(\omega_1) \rangle$ such that each $A_{\delta} = \{\alpha_{\delta, n} : n < \omega \}$ where $\alpha_{\delta, n}$'s are increasing cofinal in $\delta$ and for every $A \in [\omega_1]^{\aleph_1}$ and $b < a$, there exists some $\delta$ such that 

$$ \liminf_n \frac{|\{k < n: \alpha_{\delta, k} \in A\}|}{n} \geq b$$

\end{defn}

\begin{thm}
\label{infminus}
Let $0 < a \leq 1$ and suppose for every $b < a$, $\club^{\inf \geq b}$ holds. Then $\club^{\inf > a -}$ holds.
\end{thm}

We need two lemmas. 

\begin{lem}
\label{split}
Suppose $\club^{\inf \geq a}_S$ holds. Then there exists a partition $\langle S_i : i < \omega_1 \rangle$ of $S$ into stationary sets such that for every $i < \omega_1$, $\club^{\inf \geq a}_{S_i}$ holds.
\end{lem}

Proof of Lemma \ref{split}: Fix a witness $\bar{A} = \langle A_{\delta}: \delta \in S \rangle$ for $\club^{\inf \geq a}_S$ where each $A_{\delta} = \{\alpha_{\delta, n}: n < \omega\}$ and $\alpha_{\delta, n}$'s are increasing cofinal in $\delta$. Note that if $a \in (0.5, 1]$, this is easy - Choose $\langle X_i : i < \omega_1 \rangle$ where $X_i$'s are pairwise disjoint unbounded subsets of $\omega_1$ and let 

$$S_i = \{\delta \in S: \liminf_{n} \frac{|\{k < n: \alpha_{\delta, k} \in X_i\}|}{n} \geq a\}$$

Since $a > 0.5$, $S_i$'s are pairwise disjoint and for every $Y \in [X_i]^{\aleph_1}$, there are stationary many $\delta \in S_i$ such that 

$$\liminf_{n} \frac{|\{k < n: \alpha_{\delta, k} \in Y\}|}{n} \geq a$$

Fix $i < \omega_1$ and let $\{\alpha_{\xi}: \xi < \omega_1\}$ list $X_i$ in increasing order. Choose a club $E \subseteq \omega_1$ such that for every $\delta \in E$, $\sup_{\xi < \delta} \alpha_{\xi} = \delta$. Define $\bar{C} = \langle C_{\delta}: \delta \in S_i \rangle$ as follows. If $\delta \in E \cap S_i$, put $C_{\delta} = \{\xi: \alpha_{\xi} \in A_{\delta}\}$, otherwise choose $C_{\delta}$ arbitrarily. It is clear that $\bar{C}$ witnesses $\club^{\inf \geq a}_{S_i}$. \\

In the general case, $S_i$'s may not be pairwise disjoint but for any $F \in [\omega_1]^{K}$, where $Ka > 1$, we have $\bigcap_{i \in F} S_i = \emptyset$. For $Y \subseteq \omega_1$, let $S(Y)$ be the set of $\delta \in S$ such that 

$$\liminf_{n} \frac{|\{k < n: \alpha_{\delta, k} \in Y\}|}{n} \geq a$$

\begin{claim}
\label{wclaim}
There exists $\langle Y_i: i \in W \rangle$ such that $W \in [\omega_1]^{\aleph_1}$, each $Y_i \in [X_i]^{\aleph_1}$ and for every $i \in W$ and $Z \in [Y_i]^{\aleph_1}$, $S(Z) \setminus \bigcup_{j \in W \cap i} S(Y_j)$ is stationary.
\end{claim}

Proof of Claim \ref{wclaim}: Let $\cal{F}$ be the set of $\bar{Y} = \langle Y_i: i \in W \rangle$ where $W \in [\omega_1]^{\aleph_1}$ and each $Y_i \in [X_i]^{\aleph_1}$. For $\bar{Y} = \langle Y_i: i \in W \rangle \in \cal{F}$, let $n(\bar{Y})$ be the least $n$ such that for every $F \in [W]^n$, $\bigcap_{i \in F} S(Y_i)$ is non-stationary - So $2 \leq n(\bar{Y}) \leq K$. Let $N = \min \{n(\bar{Y}): \bar{Y} \in \cal{F}\}$ and fix $\bar{Y} = \langle Y_i: i \in W \rangle$ with $n(\bar{Y}) = N$. It suffices to show that for every $i_{\star} \in W$, there exists $j \in W$ such that $j > i_{\star}$ and for every $Z \in [Y_j]^{\aleph_1}$, $S(Z) \setminus \bigcup \{S(Y_i): i \leq i_{\star}, i \in W\}$ is stationary. Towards a contradiction, suppose this fails for some $i_{\star} \in W$. Let $W' = W \setminus (i_{\star} + 1)$. For each $j \in W'$, choose $Z_j \in [Y_j]^{\aleph_1}$ such that $S(Z_j) \setminus \bigcup \{S(Y_i): i \leq i_{\star}, i \in W\}$ is non-stationary. Let $\bar{Z} = \langle Z_j: j \in W'\rangle$. Then $n(\bar{Z}) \geq N$, so we can find $F \in [W']^{N-1}$ and  such that $\bigcap_{j \in F} S(Z_j)$ is stationary. It follows that there exists $i \in W$ such that $i \leq i_{\star}$ and $\bigcap_{j \in F} S(Z_j) \cap S(Y_i) $ is stationary. Hence $\bigcap_{j \in F \cup \{i_{\star}\}} S(Y_j)$ is also stationary: Contradiction. \qed \\

Let $\langle Y_i : i \in W \rangle$ be as in Claim \ref{wclaim}. For $i \in W$, let $T_i = S(Y_i) \setminus \bigcup_{j \in W \cap i} S(Y_j)$. Then each $T_i$ is stationary and for every $Z \in [Y_i]^{\aleph_1}$, there are stationary many $\delta \in T_i$ such that 

$$\liminf_{n} \frac{|\{k < n: \alpha_{\delta, k} \in Z\}|}{n} \geq a$$

We can now proceed as before to get a $\club^{\inf \geq a}_{T_i}$ witnessing sequence from $\langle A_{\delta}: \delta \in T_i \rangle$. This completes the proof of Lemma \ref{split}. \qed \\

\begin{lem}
\label{s1s2}
Suppose $\club^{\inf \geq a}_S$ holds and $S = S_1 \cup S_2$. Then one of $\club^{\inf \geq a}_{S_1}$, $\club^{\inf \geq a}_{S_2}$ holds.
\end{lem}

Proof of Lemma \ref{s1s2}: Fix a witness $\bar{A} = \langle A_{\delta}: \delta \in S \rangle$ for $\club^{\inf \geq a}_S$ where each $A_{\delta} = \{\alpha_{\delta, n}: n < \omega\}$ and $\alpha_{\delta, n}$'s are increasing cofinal in $\delta$.  Suppose $\club^{\inf \geq a}_{S_1}$ fails and choose $A \in [\omega_1]^{\aleph_1}$ such that for every $\delta \in S_1$

$$\liminf_{n} \frac{|\{k < n: \alpha_{\delta, k} \in A\}|}{n} < a$$

Since $\bar{A}$ is $\club^{\inf \geq a}_S$ witnessing sequence, it follows that for every $B \in [A]^{\aleph_1}$, there are stationary many $\delta \in S_2$ such that 

$$\liminf_{n} \frac{|\{k < n: \alpha_{\delta, k} \in B\}|}{n} \geq a$$

Now we can construct a $\club^{\inf \geq a}_{S_2}$ witnessing sequence as above. \qed \\

Proof of Theorem \ref{infminus}: Let $\langle a_n: n < \omega \rangle$ be an increasing sequence with $\lim_n a_n = a$. For each $n$, using Lemma \ref{split}, choose a sequence $\langle S_{n, i} : i < \omega_1\rangle$ of pairwise disjoint stationary sets such that $\club^{\inf \geq a_n}_{S_{n,i}}$ holds. For $m < n < \omega$, define $W_{m, n} = \{i < \omega_1: \club^{\inf \geq a_n}_{S_{m, i}} \text{ holds}\}$. \\

First suppose that for some $m < \omega$, there are infinitely many $n > m$ such that $W_{m, n}$ is infinite. Let $\langle n(k): k < \omega \rangle$ list such $n$'s in increasing order. Inductively choose $i(k) \in W_{m, n(k)}$ such that $i(k)$'s are pairwise distinct and $\club^{\inf \geq a_{n(k)}}_{S_{m, i(k)}}$ holds. Since $\langle S_{m, i(k)}: k < \omega \rangle$ consists of pairwise disjoint sets, the result follows. \\

So we can assume that there is no such $m$. Inductively choose a strictly increasing sequence $\langle m(k): k < \omega \rangle$ such that for every $n \geq m(k+1)$, $W_{m(k), n}$ is finite. Let $W = \bigcup \{W_{m(j), m(k)}: j < k < \omega\}$ and choose $i > \sup(W)$. Put $T_k = S_{m(k), i} \setminus \bigcup_{l < k} S_{m(l), i}$ and $T'_k = S_{m(k), i} \setminus T_k$. Then $T_k$'s are pairwise disjoint, $S_{m(k), i} = T_k \cup T'_k$ and by our choice of $i$, $\club^{\inf \geq a_{m(k)}}_{T'_k}$ does not hold. Hence, by Lemma \ref{s1s2}, $\club^{\inf \geq a_{m(k)}}_{T_k}$ must hold and we are done. \qed \\

Proof of Theorem \ref{infthms}(2): Fix $0 < a < 1$. We indicate the essential changes in the proof of Theorem \ref{infthms}(1) to get a model of $\club^{\inf \geq a} \wedge (\forall b \in (a, 1])\neg \club^{\inf \geq b}$. Define a modified countable join as follows. In Definition \ref{cp}, replace Clause (3)(ii)(b) by (b$_{\star}$) and Clause (4)(ii)(b) by (b$_{\star \star}$) below.

\begin{itemize}

\item[(b$_{\star}$)] $|\{i \in [2, j_1): (\exists k \in [m, n))(\frak{d}'_{i, k} \neq \frak{d}_{i, k}) \}| \leq j_1(1 - a)$ for every $2 < j_1 \leq j$.

\item[(b$_{\star \star}$)] $|\{i \in [N, j_1): (\exists k \in [m, n))(\frak{d}'_{i, k} \neq \frak{d}_{i, k}) \}| \leq (j_1 - N)(1 - a)$ for every $N < j_1 \leq j$.

\end{itemize}

Note that this gives rise to a transitive $\Sigma'_p$ there. Lemma \ref{limit} gets modified to the following.

\begin{lem}
\label{limit1}
Let $(\CR'_p, \Sigma'_p)$ be as in Definition \ref{cp} with (b$_{\star}$) in place of Clause (3)(ii)(b) and (b$_{\star \star}$) in place of Clause (4)(ii)(b). Let $(\CR, \Sigma)$ be an $\aleph_1$-CP such that $\CR'_p = \{\frak{c} \in \CR: \dom(\frak{c}) \subseteq \delta\}$ and $\Sigma'_p = \Sigma \res \CR'_p$. Let $\bb{Q} = \bb{Q}_{\CR, \Sigma}$, $p = \{\frak{c}^{\star}_1 = \oplus_{i \geq 1} \bar{\frak{d}}_i \}$ and $p_i = \{\frak{d}_{i, k} : k < n\}$. Then 

$$p \Vdash_{\bb{Q}} \liminf_{j} \frac{|\{i < j: p_i \in G_{\bb{Q}}\}|}{j} \geq a$$

\end{lem}

Next, Lemma \ref{clubstar} gets replaced by the following.

\begin{lem}
\label{clubb}
For every $b \in (a, 1]$, $V^{\bb{P}_{\omega_2}} \models \neg \club^{\inf \geq b}$.
\end{lem}

Proof of Lemma \ref{clubb}: Fix $b' \in (a, 1]$. Towards a contradiction, suppose $p_0 \in \bb{P}_{\omega_2}$, $\langle \name{A}_{\delta} = \{\name{\alpha}_{\delta, n} : n < \omega\} : \delta \in \Lim(\omega_1) \rangle \in V^{\bb{P}_{\omega_2}}$ are such that $p_0 \Vdash ``(\forall \delta \in \Lim(\omega_1))(\{\name{\alpha}_{\delta, n} : n < \omega\}$ is increasing cofinal in $\delta)$ and $\langle \name{A}_{\delta} : \delta < \omega_1 \rangle$ is a $\club^{\inf \geq b'}$ witnessing sequence". Since $\bb{P}_{\omega_2}$ satisfies ccc, we can find $\gamma < \omega_2$ such that $p_0 \in \bb{P}_{\gamma}$ and each $\name{\alpha}_{\delta, n}$ is a $\bb{P}_{\gamma}$-name. Fix $b \in (a, b')$. \\

Let $\name{X} = \{\alpha < \omega_1: \name{f}_{\bb{Q}_{\gamma}} = 1\}$. Then $\name{X} \in V^{\bb{P}_{\gamma + 1}}$ and $ V^{\bb{P}_{\gamma + 1}} \models \name{X} \in [\omega_1]^{\aleph_1}$. So there exist $p_1 \in \bb{P}_{\gamma}$, $q \in \bb{Q}_{\gamma}$, $\delta \in \Lim(\omega_1)$ and $n_0 < \omega$ such that $p_1 \geq p_0$ and $(p_1, q) \Vdash_{\bb{P}_{\gamma + 1}} (\forall j \geq n_0)(|\{i < j: \name{\alpha}_{\delta, i} \in \name{X}\}| \geq jb)$. We must have that $\dom(q) \cap \delta$ is unbounded in $\delta$ otherwise we can easily extend $(p_1, q)$ to get a contradiction. By possibly extending $q$, by Definition \ref{crdefn}(D)(ii), we can assume that $q = \{\frak{c}_k : k < K_{\star}\}$ where $\sup(\dom(\frak{c}_k)) < \inf(\dom(\frak{c}_{k+1}))$ for every $k < K_{\star}-1$ and for some $K < K_{\star}$, $\dom(\frak{c}_{K})$ is an unbounded subset of $\delta$. Let $S_{\gamma}$ and $\langle \frak{c}_{\gamma, \delta}: \delta \in S_{\gamma} \rangle$ witness that $(\CR_{\gamma}, \Sigma_{\gamma})$ is a thin $\aleph_1$-CP. By Claim \ref{tailagree}, we can further assume that $\frak{c}_K = \frak{c}'_{K'}$ for some $\langle \frak{c}'_n : n \leq K' \rangle \in \Sigma(\frak{c}_{\gamma, \delta})$. \\

Let $m < n < \omega$ and $\bar{\frak{d}}_i = \langle \frak{d}_{i, k} : k < n \rangle$ for $i \geq 1$ be as in Definition \ref{cp} and $\frak{c}_{\gamma, \delta} = \oplus_{i \geq 1} \bar{\frak{d}}_i$. Then as $\langle \frak{c}'_n : n \leq K' \rangle \in \Sigma(\frak{c}_{\gamma, \delta})$, we can find $N \geq 1$ a power of $2$ such that $\frak{c}_K = \frak{c}'_{K'} = \frak{c}^{\star}_N$ in the notation of Definition \ref{cp}. \\

Choose $p_2 \in \bb{P}_{\gamma}$, $p_2 \geq p_1$, $n_{\star} > n_0$ a power of $2$ and $\alpha_{n_{\star}} > \min(\dom(\frak{c}^{\star}_N))$ such that $p_2 \Vdash_{\bb{P}_{\gamma}} \name{\alpha}_{\delta, n_{\star}} = \alpha_{n_{\star}}$. Put $c = (a+b)/2$. Let $n_{\star \star} > n_{\star}$ be a power of $2$ such that $n_{\star}/n_{\star \star} < (b - c)/(1 - c)$. Choose $p_3 \geq p_2$ and $\langle \alpha_{n}: n \in [n_{\star}, n_{\star \star}) \rangle$ such that for every $n \in [n_{\star}, n_{\star \star})$, $p_3 \Vdash_{\bb{P}_{\gamma}} \name{\alpha}_{\delta, n} = \alpha_n$. Let $F = \{\alpha_n \notin \dom(q): n \in [n_{\star}, n_{\star \star})\}$. Let $q' = q \cup  \bigcup_{\alpha \in F} \{\frak{d}_{\alpha}\}$ where $\dom(\frak{d}_{\alpha}) = \{\alpha\}$ and $f_{\frak{d}_{\alpha}}(\alpha) = 0$. If $F$ is empty, put $q' = q$. \\

Now it is possible to choose $\bar{\frak{g}} \in \Sigma(\frak{c}^{\star}_N)$ such that letting $q'' = (q' \setminus \{\frak{c}^{\star}_N\}) \cup \bar{\frak{g}}$ forces $\{n \in [n_{\star}, n_{\star \star}): \alpha_n \notin \name{X}\} \geq (1 - c)(n_{\star \star} - n_{\star})$ - We leave the details of this to the reader. This means that $(p_3, q'')$ forces that $|\{i < n_{\star \star}: \name{\alpha}_{\delta, i} \in \name{X}\}| \leq n_{\star} + c(n_{\star \star} - n_{\star}) < b n_{\star \star}$ which is a contradiction. \qed \\

Now the remainder of the proof is exactly the same except for the fact that at the end of the proof of $\club^{\inf \geq a}$, we use Lemma \ref{limit1} in place of Lemma \ref{limit}. \qed \\

Proof of Theorem \ref{infthms}(3): Let $\langle a_k: k \geq 1 \rangle$ be an increasing sequence with limit $a$. Proceed as in the proof of Theorem \ref{infthms}(2) with the following modification for countable joins. In Definition \ref{cp}, replace Clause (3)(ii)(b) by (b$^{\star}$) and Clause (4)(ii)(b) by (b$^{\star \star}$) below.

\begin{itemize}

\item[(b$^{\star}$)] $|\{i \in [2, j_1): (\exists k \in [m, n))(\frak{d}'_{i, k} \neq \frak{d}_{i, k}) \}| \leq j_1(1 - a_{k_{\star}})$ for every $2 < j_1 \leq j$.

\item[(b$^{\star \star}$)] $|\{i \in [N, j_1): (\exists k \in [m, n))(\frak{d}'_{i, k} \neq \frak{d}_{i, k}) \}| \leq (j_1 - N)(1 - a_{k_{\star}})$ for every $N < j_1 \leq j$.

\end{itemize}

The rest of the proof is similar to that of Theorem \ref{infthms}(2). We leave the details to the reader. \qed

\section{On $\club^{\sup}$}

\begin{defn}
For $a \in (0, 1]$ and $S \subseteq \Lim(\omega_1)$ stationary, the principle $\club^{\sup \geq a}_S$ says the following: There exists $\bar{A} = \langle A_{\delta} : \delta \in S \rangle$ such that 

\begin{itemize}

\item[(a)] each $A_{\delta} = \{\alpha_{\delta, n} : n < \omega \}$ and $\alpha_{\delta, n}$'s are increasing cofinal in $\delta$ and

\item[(b)] for every $A \in [\omega_1]^{\aleph_1}$, there exists $\delta \in S$ such that $$ \limsup_n \frac{|\{k < n: \alpha_{\delta, k} \in A\}|}{n} \geq a$$

\end{itemize}

As usual, if $S = \Lim(\omega_1)$, we just write $\club^{\sup \geq a}$.

\end{defn}

The following remark describes the situation in the Cohen and the random reals models.

\begin{rem}
\label{remark}
(1) Suppose $V \models \club$ and let $\bb{P}$ be the forcing for adding $\aleph_2$ Cohen reals. Then $V^{\bb{P}} \models \club^{\sup \geq 1} \wedge (\forall a > 0) \neg \club^{\inf \geq a}$. Moreover, the following fails in $V^{\bb{P}}$: There exists $\bar{A} = \langle A_{\delta}: \delta \in \Lim(\omega_1) \rangle$ where each $A_{\delta} = \{\alpha_{\delta, n} : n < \omega \}$ and $\alpha_{\delta, n}$'s are increasing cofinal in $\delta$ such that for every $A \in [\omega_1]^{\aleph_1}$ and $\e > 0$, there exists some $\delta$ such that $$ \liminf_n \frac{|\{k < n: \alpha_{\delta, k} \in A\}|}{n} \geq \e$$

(2) Suppose $V \models \club$ and let $\bb{P}$ be the forcing for adding $\aleph_2$ random reals. Then $V^{\bb{P}} \models (\forall a > 0) \neg \club^{\sup \geq a}$. Moreover, the following holds in $V^{\bb{P}}$: There exists $\bar{A} = \langle A_{\delta}: \delta \in \Lim(\omega_1) \rangle$ where each $A_{\delta} = \{\alpha_{\delta, n} : n < \omega \}$ and $\alpha_{\delta, n}$'s are increasing cofinal in $\delta$ such that for every $A \in [\omega_1]^{\aleph_1}$, there exists $\e > 0$ and $\delta$ such that for every sufficiently large $n$, $$\frac{|\{k < n: \alpha_{\delta, k} \in A\}|}{n} \geq \e$$ 

\end{rem}

We now prove Theorem \ref{supequiv} - For all $a, b \in (0, 1)$, $\club^{\sup \geq a}_S$ is equivalent to $\club^{\sup \geq b}_S$. For this, it is clearly enough to show the following.

\begin{lem}
\label{square}
Let $a \in (0, 1)$ and $a \leq b < \sqrt{a}$. Then $\club^{\sup \geq a}_S$ implies $\club^{\sup \geq b}_S$.
\end{lem}

Proof of Lemma \ref{square}: Let $\bar{A} = \langle A_{\delta}: \delta \in S\rangle$ witness $\club^{\sup \geq a}_S$. We can assume that $\bar{A}$ is not a $\club^{\sup \geq b}_S$ witnessing sequence. Choose $A \in [\omega_1]^{\aleph_1}$ such that for every $\delta \in S$, for every large enough $\alpha < \delta$

$$\frac{|A \cap A_{\delta} \cap \alpha|}{|A_{\delta} \cap \alpha|} < b$$

Let $S'$ be the set of $\delta \in S$ such that

$$\limsup_{\alpha \to \delta} \frac{|A \cap A_{\delta} \cap \alpha|}{|A_{\delta} \cap \alpha|} \geq a$$

Then $S'$ is stationary. For $\delta \in S'$, define $B_{\delta} = A_{\delta} \cap A$. 

\begin{claim} 
\label{case1}
For every $B \in [A]^{\aleph_1}$ there are stationary many $\delta \in S'$ such that
$$\limsup_{\alpha \to \delta} \frac{|B \cap B_{\delta} \cap \alpha|}{|B_{\delta} \cap \alpha|} \geq b $$

\end{claim}

Proof of Claim \ref{case1}: Suppose not. Choose $B \in [A]^{\aleph_1}$ and $W \subseteq S'$ non stationary such that for every $\delta \in S' \setminus W$, for every large enough $\alpha < \delta$, we have 

$$\frac{|B \cap B_{\delta} \cap \alpha|}{|B_{\delta} \cap \alpha|} < b $$

Since $B \subseteq A$, we can choose $\delta \in S' \setminus W$ such that 

$$\limsup_{\alpha \to \delta} \frac{|B \cap A_{\delta} \cap \alpha|}{|A_{\delta} \cap \alpha|} \geq a $$

Now for every large enough $\alpha < \delta$, we have

$$\left( \frac{|B \cap B_{\delta} \cap \alpha|}{|B_{\delta} \cap \alpha|} \right) \left( \frac{|A \cap A_{\delta} \cap \alpha|}{|A_{\delta} \cap \alpha|} \right) < b^2 $$

Since $B \cap B_{\delta} = B \cap A_{\delta}$ and $B_{\delta} \cap \alpha = A \cap A_{\delta} \cap \alpha$, we get

$$\frac{B \cap A_{\delta} \cap \alpha}{A_{\delta} \cap \alpha} < b^2 < a$$

which is impossible. \\

Let $\{\alpha_i : i < \omega_1\}$ list $A$ in increasing order. Let $E \subseteq \omega_1$ be a club such that for every $i \in E$, $\sup_{j < i} \alpha_j = i$. Define $\bar{C} = \langle C_{\delta}: \delta \in S \rangle$ as follows. If $\delta \in E \cap S'$, then $C_{\delta} = \{j < \delta: \alpha_j \in B_{\delta}\}$. Otherwise, choose $C_{\delta}$ to be an arbitrary unbounded subset of $\delta$ of order type $\omega$. It is easy to check that $\bar{C}$ witnesses $\club^{\sup \geq b}_{S}$. \qed

\section{$\neg \club^{\sup \geq 1}$ and $\club^{\sup > 1-}$}

\begin{defn}
The principle $\club^{\sup > 1 -}$ says the following: There exists $\bar{A} = \langle A_{\delta} : \delta \in \Lim(\omega_1) \rangle$ such that 

\begin{itemize}

\item[(a)] each $A_{\delta} = \{\alpha_{\delta, n} : n < \omega \}$ and $\alpha_{\delta, n}$'s are increasing cofinal in $\delta$ and

\item[(b)] for every $A \in [\omega_1]^{\aleph_1}$ and $\e > 0$, there exists some $\delta$ such that $$ \limsup_n \frac{|\{k < n: \alpha_{\delta, k} \in A\}|}{n} \geq 1 - \e$$

\end{itemize}

\end{defn}

To prove Theorem \ref{supsep}, it is enough to show that

\begin{thm}
\label{supminus}
$\neg \club^{\sup \geq 1} \wedge \club^{\sup > 1-}$ is consistent.
 
\end{thm}

\begin{defn}
\label{forcesup}
Suppose $\bar{A} = \langle A_{\delta} : \delta \in \Lim(\omega_1) \rangle$ satisfies: For every $\delta$, $A_{\delta} = \{\alpha_{\delta, n}: n < \omega\}$ where $\alpha_{\delta, n}$'s are increasing and cofinal in $\delta$. Define $\bb{Q} = \bb{Q}_{\bar{A}}$ as follows: $p \in \bb{Q}$ iff $p = (f_p, u_p, \bar{\e}_p)$ where 

\begin{itemize}

\item[(i)] $f_p$ is a finite partial function from $\omega_1$ to $\{0, 1\}$,

\item[(ii)] $u_p$ is a finite subset of $\Lim(\omega_1)$ and

\item[(iii)] $\bar{\e}_p = \langle \e_{p, \delta}: \delta \in u_p \rangle$ where each $\e_{p, \delta}$ is a positive rational $< 1$.

\end{itemize}

For $p, q \in \bb{Q}$ define $p \leq q$ iff 

\begin{itemize}

\item[(a)] $f_p \subseteq f_q$,

\item[(b)] $u_p \subseteq u_q$,

\item[(c)] $\bar{\e}_{p} = \bar{\e}_q \res u_p$ and

\item[(d)] for every $\delta \in u_p$, letting $W = \{n < \omega: \alpha_{\delta, n} \in \dom(f_q) \setminus \dom(f_p)\}$, for every $N < \omega$ either $W \cap [0, N) = \emptyset$ or

$$ \frac{|\{n \in W \cap [0, N): f_q(\alpha_{\delta, n}) = 1\}|}{|W \cap [0, N)|} \leq 1 - \e_{p, \delta}$$

\end{itemize}

\end{defn}

\begin{claim}
\label{supccc}
Let $\bar{A}$ and $\bb{Q} = \bb{Q}_{\bar{A}, a}$ be as in Definition \ref{forcesup}. Then $\bb{Q}$ has $\aleph_1$ as a precaliber.

\end{claim}

Proof of Claim \ref{supccc}: Suppose $\{p_i = (f_i, u_i, \bar{\e}_i): i < \omega_1\} \subseteq \bb{Q}$. By thinning down we can assume the following.

\begin{itemize}

\item[(a)] $\langle \dom(f_i) : i < \omega_1\rangle$ is a $\Delta$-system with root $R$ and $f_i \res R$ does not depend on $i$.

\item[(b)] $\langle u_i : i < \omega_1 \rangle$ is a $\Delta$-system with root $u_{\star}$ and $\bar{\e}_i \res u_{\star}$ does not depend on $i$.

\item[(c)] For every $i < j < \omega_1$ and $\delta \in u_i$, $\dom(f_j) \cap A_{\delta} \subseteq R$.

\end{itemize}

Let $E \subseteq \omega_1$ be a club such that for every $i \in E$, for every $j < i$, $\dom(f_j) \cup u_j \subseteq i$. Choose $S \subseteq E$ stationary such that for every $i \in S$, $\dom(f_i) \cap i = R$,  $u_i \cap i = u_{\star}$ and $\bigcup \{A_{\delta} \cap i: \delta \in u_i, \delta > i\} = F$ where $F$ does not depend on $i \in S$. Note that for every infinite $X \subseteq S$ and $i \in S$, if $i > \sup(X)$, then for all but finitely many $j \in X$, $\dom(f_j) \cap A_i \subseteq R$. Let $X \in [S]^{\aleph_1}$ be such that for every increasing sequence $\langle \alpha_n: n < \omega \rangle$ in $X$, $\sup_n \alpha_n \notin X$. Define $c:[X]^2 \to \{0, 1\}$ by $c(\{i, j\}) = 1$ iff $i < j$ and $A_j \cap \dom(f_i) \subseteq R$. By Erdos-Dushnik-Miller, either there exists $Y \in [X]^{\aleph_1}$ such that $c[[Y]^2] = \{1\}$ or there exists $Y' \subseteq X$ such that $\otp(Y') = \omega+1$ and $c[[Y']^2] = \{0\}$. Since the latter is impossible, we can find $Y \in [X]^{\aleph_1}$ such that $c[[Y]^2] = \{1\}$. Hence

\begin{itemize}

\item[(d)] For every $i \neq j$ in $Y$ and $\delta \in u_j$, $\dom(f_i) \cap A_{\delta} \subseteq R$.

\end{itemize}

It follows that $\{ p_i: i \in Y \}$ is centered. \qed \\

Let $\name{f}_{\bb{Q}} = \bigcup \{f_p : p \in G_{\bb{Q}}\}$. Then $\Vdash_{\bb{Q}} \name{f}_{\bb{Q}}: \omega_1 \to \{0, 1\}$. Let $\name{X}_{\bb{Q}} = \{\alpha < \omega_1 : \name{f}_{\bb{Q}}(\alpha) = 1\}$. Then $\Vdash_{\bb{Q}} \name{X}_{\bb{Q}} \in [\omega_1]^{\aleph_1}$.

\begin{claim}
$\name{X}_{\bb{Q}}$ witnesses that $\bar{A}$ is not a $\club^{\sup \geq 1}$ witnessing sequence in $V^{\bb{Q}}$.

\end{claim}

Proof: Easy. \qed

\begin{claim}
\label{clubsuppres}
Suppose $V \models \club^{\sup > 1-}$ holds and let $\bar{C} = \langle C_{\delta} : \delta \in \Lim(\omega_1) \rangle$ be a witness where $C_{\delta} = \{\beta_{\delta, n}: n < \omega\}$ and $\beta_{\delta, n}$'s are increasing cofinal in $\delta$. Then $V^{\bb{Q}} \models \club^{\sup > 1 -}$ holds with $\bar{C}$ as witness.

\end{claim}

Proof of Claim \ref{clubsuppres}: Suppose $p \Vdash_{\bb{Q}} \name{A} \in [\omega_1]^{\aleph_1}$ and $\e > 0$. Choose $\langle (p_i, \gamma_i): i < \omega_1 \rangle$ such that $\gamma_i$'s are increasing and for each $i < \omega_1$, $p \leq p_i \Vdash_{\bb{Q}} \gamma_i \in \name{A}$. Arguing as in the proof of Claim \ref{supccc}, we can assume the following.

\begin{itemize}

\item[(a)] $\langle \dom(f_i) : i < \omega_1\rangle$ is a $\Delta$-system with root $R$, $f_i \res R = f_{\star}$ and $|\dom(f_i) \setminus R| = n_{\star}$ do not depend on $i$.

\item[(b)] If $i < j$, then $R < \dom(f_i) \setminus R < \dom(f_j) \setminus R$.

\item[(c)] $\langle u_i : i < \omega_1 \rangle$ is a $\Delta$-system with root $u_{\star}$, $\bar{\e}_i \res u_{\star} = \bar{e}_{\star}$ does not depend on $i$ and $i < j$ implies $u_i \setminus u_{\star} < u_j \setminus u_{\star}$.

\item[(c)] For every $i \neq j$ and $\delta \in u_i$, $\dom(f_j) \cap A_{\delta} \subseteq R$.

\end{itemize}

Put $X = \{\gamma_i: i < \omega_1\}$. Let $E \subseteq \omega_1$ be a club such that for every $i \in E$ and $j < i$, $\gamma_j < i$ and $u_{\star} \cup \dom(f_j) \subseteq i$. Choose $\delta \in E$ such that 

$$\limsup_n \frac{|\{k < n: \beta_{\delta, k} \in X\}|}{n} \geq 1 - \e/10$$

Let $q = (f_{\star}, u_{\star} \cup \{\delta\}, \bar{\e}_{\star} \cup \{(\delta, \e/5)\})$. It suffices to show that for any $q_1 \geq q$ and $N_0 < \omega$, there exist $r \geq q_1$ and $N_2 > N_0$ such that 

$$r \Vdash_{\bb{Q}} \frac{|\{n < N_2: \beta_{\delta, n} \in \name{A}\}|}{N_2} \geq 1 - \e$$

So fix $q_1 \geq q$ and $N_0 < \omega$. For each $n < \omega$, define

\[
  r_n =
  \begin{cases}
    p_i & \text{if } \beta_{\delta, n} = \gamma_i \\
    q & \text{if } \beta_{\delta, n} \notin X
  \end{cases}
\]

Let $W'_n = \dom(f_{r_n}) \setminus R$ and $W_n = W'_n \cap A_{\delta}$. Choose $N_1 > N_0$ such that for every $n \geq N_1$, if $\delta' \in u_{q_1} \setminus \{\delta\}$, then $W'_n \cap A_{\delta'} = \phi$. We need a lemma.

\begin{lem}
\label{count}
Suppose $0 < a_1 < a_2 < 1$ and $1 \leq K < \omega$. Then for all sufficiently large $N < \omega$, the following holds. For every $\langle W_k : k < N \rangle$ where each $W_k$ is an interval in $\omega$ such that $|W_k| \leq K$, $W_k < W_{k+1}$ and $\bigcup_{k < n} W_k = [0, M)$, there exists $F \subseteq N$ such that

\begin{itemize}

\item[(i)] $|F| \geq Na_1$ and

\item[(ii)] For every $m \leq M$, $|[0, m) \cap \bigcup_{k \in F} W_k| \leq m a_2$

\end{itemize}

\end{lem}

Proof of Lemma \ref{count}: First assume that $|W_k| =  K$ for every $k < N$ - So $M = NK$. Let $m_1 < N$ be least such that $K m_1 \geq M (1 - a_2)$. Then $F = [m_1, N)$ is as required. For the general case, for each $K' \leq K$, put $S_{K'} = \{k < N : |W_k| = K'\}$ and find a suitable $F_{K'} \subseteq S_{K'}$ for $\langle W_k : k \in S_{K'} \rangle$. Then $F = \bigcup \{F_{K'}: 1 \leq K' \leq K\}$ is as required. \qed \\

Choose $N_2 > N_1$ such that $(1 - N_1/N_2)(1 - \e/2) \geq 1 - \e$ and $|\{k \in [N_1, N_2): \beta_{\delta, k} \in X\}| \geq (1 - \e/4)(N_2 - N_1)$. Using Lemma \ref{count}, choose $F \subseteq [N_1, N_2)$ such that the following hold.

\begin{itemize}

\item[(a)] $|F| \geq (N_2 - N_1)(1 - \e/4)$.

\item[(b)] $r = (f_r, u_r, \bar{\e}_r)$ extends each condition in $\{q_1, r_n: n \in F \}$ where 

\subitem (i) $u_r = u_{q_1} \cup \bigcup_{n \in F} u_{r_n}$,

\subitem (ii) $\dom(f_r) = \dom(f_{q_1}) \cup \bigcup_{n \in F} W'_n \cup \bigcup \{W_n : n \in [N_1, N_2) \setminus F\}$ ,

\subitem (iii) $f_{q_1} \subseteq f_r$,

\subitem (iv) $f_r \res \bigcup \{W_n : n \in [N_1, N_2) \setminus F\} \equiv 0$,

\subitem (v) for every $n \in F$, $f_r \res W'_n = f_{r_n}$ and

\subitem (vi) $\bar{\e}_r = \bar{\e}_{q_1} \cup \bigcup_{n \in F} \bar{\e}_{r_n}$.

\end{itemize}

Note that $r \Vdash |\{k < N_2: \beta_{\delta, k} \in \name{A}\}| \geq (N_2 - N_1)(1 - \e/2)$. By our choice of $N_2$, it follows that 

$$r \Vdash_{\bb{Q}} \frac{|\{n < N_2: \beta_{\delta, n} \in \name{A}\}|}{N_2} \geq 1 - \e$$ \qed

Let $\eta \geq 1$ and suppose $\langle (\bb{P}_{\xi} , \bb{Q}_{\xi}, \bar{A}_{\xi}) : \xi < \eta \rangle$ satisfies the following.

\begin{itemize}

\item[(1)] $\langle (\bb{P}_{\xi} , \bb{Q}_{\xi}) : \xi < \eta \rangle$ is a finite support iteration with limit $\bb{P}_{\eta}$.

\item[(2)] $\bar{A}_{\xi} \in V^{\bb{P}_{\xi}}$ and $ \Vdash_{\bb{P}_{\xi}} ``\bar{A}_{\xi} = \langle A_{\xi, \delta} : \delta \in \Lim(\omega_1) \rangle$, $A_{\xi, \delta} = \{\alpha_{\xi, \delta, n} : n < \omega\}$ where $\alpha_{\xi, \delta, n}$'s are increasing cofinal in $\delta$".

\item[(3)] $V^{\bb{P}_{\xi}} \models \bb{Q}_{\xi} = \bb{Q}_{\bar{A}_{\xi}}$.

\end{itemize}

Note that $\bb{P}_{\eta}$ is ccc.

\begin{claim}
\label{clubsuppres1}
Suppose $V \models \club^{\sup > 1-}$ holds and let $\bar{C} = \langle C_{\delta} : \delta \in \Lim(\omega_1) \rangle$ be a witness where $C_{\delta} = \{\beta_{\delta, n}: n < \omega\}$ and $\beta_{\delta, n}$'s are increasing cofinal in $\delta$. Then $V^{\bb{P}_{\eta}} \models \club^{\sup \geq 1-}$ via the same witness.
\end{claim}

Proof of Claim \ref{clubsuppres1}: By induction on $\eta$. If $\eta$ is a successor or $\cf(\eta) > \aleph_1$, this follows from Claim \ref{clubsuppres}. \\

Suppose $\cf(\eta) = \aleph_0$. Let $\langle \eta(n) : n < \omega \rangle$ be increasing cofinal in $\eta$. Suppose $p \Vdash_{\bb{P}_{\eta}} \name{X} \in [\omega_1]^{\aleph_1}$. Choose $n_{\star} < \omega$ such that $p \in \bb{P}_{\eta(n_{\star})}$ For each $n < \omega$, let $\name{X}_n = \{\alpha < \omega_1: (\exists p \in G_{\bb{P}_{\eta(n)}})(p \Vdash_{\bb{P}_{\eta}} \alpha \in \name{X}) \}$ - So $\name{X}_n \in V^{\bb{P}_{\eta(n)}}$ and $\Vdash_{\bb{P}_{\eta}} \name{X}_n \subseteq \name{X}$. Then for some $n \in [n_{\star}, \omega)$, $p \Vdash_{\bb{P}_{\eta(n)}} \name{X}_n \in [\omega_1]^{\aleph_1}$. Now apply the inductive hypothesis. \\

Next suppose $\cf(\eta) = \aleph_1$, $\e > 0$ and $p \Vdash_{\bb{P}_{\eta}} \name{X} \in [\omega_1]^{\aleph_1}$. Choose $\langle (p_i, \gamma_i): i < \omega_1 \rangle$ such that the following hold.

\begin{itemize}

\item[(a)] $\gamma_i$'s are increasing.

\item[(b)] $p_i \in \bb{P}_{\eta}$, $p_i \geq p$ and $p_i \Vdash_{\bb{P}_{\eta}} \gamma_i \in \name{X}$.

\item[(c)] $\langle \dom(p_i): i < \omega_1 \rangle$ is a $\Delta$-system with root $W$.

\end{itemize}

Choose $\theta < \eta$ such that $W \subseteq \theta$. Since $\bb{P}_{\theta}$ is ccc, we can find $q \in \bb{P}_{\theta}$ such that $q \geq p$ and $q \Vdash_{\bb{P}_{\theta}} ``\{i < \omega_1: p_i \res \theta \in G_{\bb{P}_{\theta}} \}$ is uncountable". Let $\name{Y} = \{\gamma_i: i < \omega_1 \wedge p_i \res \theta \in G_{\bb{P}_{\theta}}\}$. Then $\name{Y} \in V^{\bb{P}_{\theta}}$ and $q \Vdash_{\bb{P}_{\theta}} \name{Y} \in [\omega_1]^{\aleph_1}$. By the inductive hypothesis, we can find $r \in \bb{P}_{\theta}$ and $\delta \in \Lim(\omega_1)$ such that $r \geq q$ and 

$$r \Vdash_{\bb{P}_{\theta}} \limsup_n \frac{|\{k < n: \beta_{\delta, k} \in \name{Y}\}|}{n} \geq 1 - \e/2$$

Since $\langle \dom(p_i) \setminus \theta : i < \omega_1 \rangle$ is a sequence of pairwise disjoint sets, it also follows that

$$r \Vdash_{\bb{P}_{\eta}} \limsup_n \frac{|\{k < n: \beta_{\delta, k} \in \name{X}\}|}{n} \geq 1 - \e$$

\qed

Proof of Theorem \ref{supminus}: Starting with a model of $2^{\aleph_1} = \aleph_2$ and $\club^{\sup > 1-}$ construct $\langle (\bb{P}_{\xi} , \bb{Q}_{\xi}, \bar{A}_{\xi}) : \xi < \omega_2 \rangle$ such that the following hold.

\begin{itemize}

\item[(1)] $\langle (\bb{P}_{\xi} , \bb{Q}_{\xi}) : \xi < \omega_2 \rangle$ is a finite support iteration with limit $\bb{P}_{\omega_2}$.

\item[(2)] $\bar{A}_{\xi} \in V^{\bb{P}_{\xi}}$ and $ \Vdash_{\bb{P}_{\xi}} ``\bar{A}_{\xi} = \langle A_{\xi, \delta} : \delta \in \Lim(\omega_1) \rangle$, $A_{\xi, \delta} = \{\alpha_{\xi, \delta, n} : n < \omega\}$ where $\alpha_{\xi, \delta, n}$'s are increasing cofinal in $\delta$".

\item[(3)] $V^{\bb{P}_{\xi}} \models \bb{Q}_{\xi} = \bb{Q}_{\bar{A}_{\xi}}$.

\item[(4)] For every $\eta < \omega_2$ and $\bar{A} \in V^{\bb{P}_{\eta}}$ satisfying $\Vdash_{\bb{P}_{\eta}} ``\bar{A} = \langle A_{\delta}: \delta \in \Lim(\omega_1) \rangle$ where each $A_{\delta}$ is an unbounded subset of $\delta$ of order type $\omega$", there exists $\xi \in [\eta, \omega_2)$ such that $\Vdash_{\bb{P}_{\xi}} \bar{A} = \bar{A}_{\xi}$.

\end{itemize}

To see why clause (4) can be satisfied, use $2^{\aleph_1} = \aleph_2$ and the fact that for each $\eta < \omega_2$, $\bb{P}_{\eta}$ is a ccc forcing with a dense subset of size $\aleph_1$. \qed \\

We conclude with some questions.
 
\begin{ques}

\begin{itemize}

\item[(1)] Is $\club^{\sup \geq 0.5} \wedge \neg \club^{\sup > 1 -}$ consistent? What if CH holds?

\item[(2)] Assume CH. Does $\club^{\sup \geq 0.5}$ imply $\club^{\sup \geq 1}$? Does $\club^{\sup > 1-}$ imply $\club^{\sup \geq 1}$?

\item[(3)] For $a \in (0, 1)$, is $\club^{\inf \geq a} \wedge \neg \club^{\sup \geq 1}$ consistent?

\end{itemize}

\end{ques}

\end{document}